\documentclass[11pt]{article}

\usepackage{array} 
\usepackage{amsfonts}
\usepackage{amsmath}
\usepackage{amssymb}
\usepackage{graphicx}
\usepackage{color}

\newtheorem{thm}{Theorem}[section]
\newtheorem{cor}[thm]{Corollary}
\newtheorem{lem}[thm]{Lemma}
\newtheorem{prop}[thm]{Proposition}
\newtheorem{defn}[thm]{Definition}

\numberwithin{equation}{section}

\newcommand{\dx}{\,{\rm d}x}
\newcommand{\dy}{\,{\rm d}y}

\newcommand{\dt}{\,{\rm d}t}
\newcommand{\rd}{{\rm d}}

\def\LL{\mathrm{L}} 
\def\supp{\mathrm{supp}} 

\newcommand{\A}{\mathcal{L}}

\newcommand{\AI}{\mathcal{L}^{-1}}

\newcommand{\n}{F}
\newcommand{\nl}{{F^*}}

\newcommand{\p}{{\delta_\gamma}} 

\newcommand{\ka}{\overline{\kappa}}
\newcommand{\kb}{\underline{\kappa}}

\newcommand{\K}{{\mathbb K}}

\newcommand{\RR}{\mathbb{R}}

\def\ee{\mathrm{e}} 
\def\dist{\mathrm{dist}} 
\def\dom{\mathrm{dom}} 

\def\qed{\,\unskip\kern 6pt \penalty 500
\raise -2pt\hbox{\vrule \vbox to8pt{\hrule width 6pt
\vfill\hrule}\vrule}\par}
\definecolor{darkblue}{rgb}{0.05, .05, .65}
\definecolor{darkgreen}{rgb}{0.1, .65, .1}
\definecolor{darkred}{rgb}{0.8,0,0}

\begin{document}
\title{\textbf{Fractional Nonlinear Degenerate \\ Diffusion Equations  on Bounded Domains\\ Part I. Existence, Uniqueness and Upper Bounds}\\[7mm]}

\author{{\Large Matteo Bonforte\footnote{e-mail address:~\texttt{matteo.bonforte@uam.es}} ~and~ Juan Luis V\'azquez\footnote{e-mail address:~\texttt{juanluis.vazquez@uam.es}}}\\ [4pt]
Departamento de Matem\'{a}ticas \\ [4pt] Universidad
Aut\'{o}noma de Madrid} 
\date{} 

\maketitle

\begin{abstract}
We investigate  quantitative properties of nonnegative solutions $u(t,x)\ge 0$ to the nonlinear fractional diffusion equation, $\partial_t u
+ \A \n(u)=0$ posed in a bounded domain, $x\in\Omega\subset \RR^N$\,, with appropriate homogeneous Dirichlet boundary conditions. As $\A$ we can use a quite general class of linear operators that includes the two most common versions of the fractional Laplacian $(-\Delta)^s$, $0<s<1$, in a bounded domain with zero Dirichlet boundary conditions, but it  also includes many other examples since our theory only needs some basic properties that are typical of  ``linear heat semigroups''. \normalcolor The nonlinearity $\n$ is assumed to be increasing and is allowed to be degenerate, the prototype is the power case $\n(u)=|u|^{m-1}u$, with $m>1$\,.

In this paper we propose a  suitable class of  solutions of the equation, and cover the basic theory: we prove existence, uniqueness of such solutions, and we establish upper bounds of two forms (absolute bounds and smoothing effects), as well as weighted-$\LL^1$ estimates. The class of solutions is very well suited for that work. The standard Laplacian case $s=1$ is included and the linear case $m=1$ can be recovered in the limit.

In a companion paper \cite{BV-ppr2-2}, we will complete the study with more advanced estimates, like the upper and lower boundary behaviour and Harnack inequalities, for which the results of this paper are needed.
\end{abstract}\vspace{2cm}


\newpage
\small
\tableofcontents
\normalsize
\newpage

\section{Introduction}

In this paper we address the question of  obtaining an existence and uniqueness theory together with quantitative  a priori estimates for a suitable class of weak solutions of the nonlinear  fractional diffusion equation of the form:
\begin{equation}\label{FPME.equation}
\partial_t u+\A\,\n (u)=0\,,
\end{equation}
posed in $Q=(0,\infty)\times \Omega$ where $\Omega\subset \RR^N$ is a bounded domain with smooth boundary, $N\ge 1$,  $\A$ is a linear operator representing diffusion, mainly of the fractional Laplacian type, and $F:\RR \to \RR$ is a monotone nondecreasing real function. We will refer to the equation as (NFDE).  Since the equation is posed in a bounded domain, we need boundary conditions that we assume of Dirichlet type and must be  defined in accordance to the nonlocal nature of the fractional operator.

This kind of problem has been extensively studied when $\A=-\Delta$ and $F$ is a power function, $\n(u)=|u|^{m-1}u$, and then the equation becomes the Porous Medium Equation \cite{VazBook} if $m>1$, see also \cite{DK0, DaskaBook, JLVmonats}. For $\n(u)=u$ it is simply the Heat Equation. We are interested here in treating nonlocal diffusion operators, in particular the fractional Laplacian operators. Note that, since we are working on a bounded domain,  the latter concept admits  at least two non-equivalent versions, the Restricted Fractional Laplacian (RFL) and the Spectral Fractional Laplacian (SFL),  and then many variants. The case of the SFL operator and  $F$ a power function with $m>1$  was already studied by us in  \cite{BV-PPR1}.  Starting from that precedent, we present as the contribution of the present paper a rather abstract setting where we are able to treat a quite general class of operators, in particular singular integral operators with more general kernels, and we consider nonlinearities $F$ under a list of hypotheses that include the case of power functions.
We refer to \cite{AC,BV-PPR1,BV2012,DPQRV1,DPQRV2} for the physical motivation and relevance of this nonlocal model.

A brief outline of the paper is as follows. We first introduce the precise conditions on the linear operator and the nonlinearity, and we check that the intended examples fall into this setting. The next step is to propose a class of generalized solutions of equation \eqref{FPME.equation}, called {\sl weak dual solutions}, that have been introduced in paper \cite{BV-PPR1} as the most  convenient class for these problems.

We want to establish a theory of existence and uniqueness for such equations and such solutions. It relies on an abstract existence theorem plus a set of a priori estimates. The first part is based on the abstract theory of mild solutions constructed by Crandall and Pierre \cite{CP-JFA} for $m$-accretive operators\,, that we recall in Theorem \ref{thm.CP-JFA}.

 This is followed by three  sections establishing a list of a priori estimates for nonnegative dual weak solutions of equation (NFDE) under our assumptions. The long sections \ref{sect.frist.two}, \ref{sect.smooth} and \ref{Sect.Weighted.L1} represent the core of this paper. Precisely, we prove absolute upper bounds, smoothing effects and weighted $\LL^1$ estimates.

Once this is accomplished, we have to check that both parts agree. Actually,  we have to prove that the semigroup (mild) solutions constructed in Theorem \ref{thm.CP-JFA} are indeed weak dual solutions, at least for a nice class of initial data  (which is $L^1(\Omega)$), therefore they inherit the a priori estimates valid for the dual weak solutions. The various definitions of solutions and the relations among them is treated in Section \ref{Sect.Various.Sols}.  Finally in Section   \ref{sect.exist.uniq.WDS}, we prove existence and uniqueness of minimal dual weak solutions with initial data in    a weighted $\LL^1$ space, by approximation.

We collect in the Appendix a number of technical lemmas used in the proofs. A final section of the paper contains  comments, consequences, ideas on extensions, and related works.

Let us comment a related work of the authors with Y. Sire \cite{BSV2013} where a similar problem setting is discussed, as a preliminary step to perform the asymptotic analysis of solutions as $t\to\infty$. In such paper the functional framework for existence is the variational setting $H^{-s}(\Omega)$ which leads to a class of weak solutions in Hilbert spaces in the spirit of Brezis \cite{Brezis71}. By contrast our present theory wants to include the classes of solutions in $L^1(\Omega)$ typical of the $L^1$-semigroup approach of \cite{CL71, CP-JFA}. The new class of weak dual solutions is particularly friendly for obtaining a priori estimates. This method can be extended to a wider class of nonlocal diffusion operators $\A$. Let us also point out that the method we use leads naturally to estimates  in the weighted space $L^1_{\p}(\Omega)$, which allows to include the $H^{-s}(\Omega)$ solutions (when $u\ge 0$, as we assume in this paper).

\medskip

It is maybe worth making a comment about the novelty of the paper. On one hand, we exploit the functional properties of the linear operator $\cal L$ as much as possible. A key ingredient is the knowledge of good estimates for the Green function, that hold for the standard examples of fractional operator, and also for many variants. This information is essential in formulating the definition of weak dual solution following the trend set in our previous works \cite{BV-PPR1, BSV2013}; it allows us here to present a theory that adapts to a wide class of operators. The core of the paper is devoted to use that setting to derive very precise a priori estimates. The third goal of the article is an abstract existence and uniqueness theory that is shown to be compatible with the previous items. Summing up, the functional approach allows us to avoid more standard  and delicate methods; in this way it is possible  to treat a quite wide class of linear operators $\cal L$ combined with nonlinearities $F$.    Such generality should be explored in future works.

\section{First assumptions and notations}

\noindent $\bullet$ {\bf On $F$.} We will always consider $\n:\RR\to\RR$ to be a continuous and non-decreasing function, with the  normalization $F(0)=0$. Moreover, it satisfies the condition:
\begin{enumerate}
\item[(N1)] $\n\in C^1(\RR\setminus\{0\})$ and $\n/\n'\in {\rm Lip}(\RR)$ and there exists $\mu_0,\mu_1>0$ such that
\[
1-\mu_1\le \left(\frac{\n}{\n'}\right)'\le 1-\mu_0\,,
\]
where $\n/\n'$ is understood to vanish if $\n(r)=\n'(r)=0$ or $r=0$\,.
\end{enumerate}
The main example will be $\n(u)=|u|^{m-1}u$, with $m>1$\,, which corresponds to the nonlocal porous medium equation and has been studied in \cite{BV-PPR1}. Then $\mu_0=\mu_1=(m-1)/m$.  A simple variant is the combination of two powers, so that one of them gives the behaviour near $u=0$, the other one the behaviour near $u=\infty$.

\medskip

\noindent $\bullet$ {\bf On $\A$.} The linear operator $\A: \dom(A)\subseteq\LL^1(\Omega)\to\LL^1(\Omega)$ is assumed to be densely defined  and sub-Markovian, more precisely satisfying $(A1)$ and $(A2)$ below:
\begin{enumerate}
\item[(A1)] $\A$ is $m$-accretive on $\LL^1(\Omega)$,
\item[(A2)] If $0\le f\le 1$ then $0\le \ee^{-t\A}f\le 1$\,,
    or equivalently,
\item[(A2')] If $\beta$ is a maximal monotone graph in $\RR\times\RR$ with $0\in \beta(0)$, $u\in \dom(\A)$\,, $\A u\in \LL^p(\Omega)$\,, $1\le p\le\infty$\,, $v\in \LL^{p/(p-1)}(\Omega)$\,, $v(x)\in \beta(u(x))$ a.e\,, then
    \[
    \int_\Omega v(x)\A u(x)\dx\ge 0
    \]
\end{enumerate}

These assumptions are taken from Crandall and Pierre's paper \cite{CP-JFA} where semigroup (mild) solutions are constructed for the abstract equation $u_t+\A\n(u)=0$. The proper interpretation of $\A\n$ has been given there. Actually, (A1) and (A2) imply that $\dom(\A)$ is dense,  cf. \cite{CP-JFA}. Another interesting result of their paper is represented by the monotonicity estimates. We summarize next some results of \cite{CP-JFA} that we will use in the rest of this paper.

\begin{thm}[Crandall-Pierre, \cite{CP-JFA}]\label{thm.CP-JFA}
Let $\A$ satisfy $(A1)$ and $(A2)$ and let $\n$ satisfy $(N1)$. Then for all nonnegative $u_0\in \LL^1(\Omega)$\,, there exists a unique mild solution $u$ to equation \eqref{FPME.equation}\,, and the function
\begin{equation}\label{CP.Monotonicity}
t\mapsto t^{\frac{1}{\mu_0}}\,\n(u(t,x))\qquad\mbox{is nondecreasing in $t>0$ for a.e. }x\in \Omega\,.
\end{equation}
Moreover, the semigroup is contractive on $\LL^1(\Omega)$ and $u\in C([0,\infty): \LL^1(\Omega))$\,.
\end{thm}
We notice that \eqref{CP.Monotonicity} is a weak formulation of the monotonicity inequality:
\begin{equation}\label{CP.Monotonicity.2}
\partial_t u \ge -\frac{1}{\mu_0\, t}\frac{\n(u)}{\n'(u)}
\end{equation}
(where $\n(u)/\n'(u)$ is understood to vanish if $\n(u)=\n'(u)=0$ or if $u=0$). Moreover, condition $(N2)$ implies that $\n(u)/\n'(u)\le (1-\mu_0)u$ so that the above inequality becomes:
\begin{equation}\label{CP.Monotonicity.3}
\partial_t u \ge -\frac{1-\mu_0}{\mu_0}\,\frac{u}{t}
\end{equation}
or equivalently the function $t\mapsto t^{\frac{1-\mu_0}{\mu_0}}\,u(t,x)$ is nondecreasing in $t>0$ for a.e. $x\in \Omega$\,.

We shall remark that when $\n(r)=|r|^{m-1}r$\,, we can take $\mu_0=(m-1)/m<1$ to recover the well known monotonicity estimates by B\'enilan and Crandall \cite{BCr}:
\[
\partial_t u \ge -\frac{u}{(m-1)t}
\]

\noindent {\bf Remark.} The solutions constructed in \cite{CP-JFA} are mild solutions, and  it is not easy to prove properties like boundedness, positivity,  or to establish boundary estimates for that concept of solution, that arises  in the semigroup setting.  One of the goals of the present  paper and of the forthcoming one, \cite{BV-ppr2-2}, is to show that such solutions are indeed bounded, and positive after a waiting time. Moreover,  we will   prove that close to $\partial\Omega$ they behave as  $\n^{-1}$ of the distance to the boundary. Indeed, under some further assumption on $\A$ (but not on $\n$), we are going to construct a larger class of solutions, that we will call \textit{weak dual solutions} which satisfy all the above mentioned ``good'' properties.

\subsection{Additional assumptions on $\A$}

In other to construct our version of the theory and prove the desired quantitative properties, we need to be more specific about operator $\A$. Besides satisfying (A1) and (A2), we will assume that it has a left-inverse $\AI: \LL^1(\Omega)\to \LL^1(\Omega)$  with a kernel $\K$ such that
\[
\AI[f](x)=\int_\Omega \K(x,y)f(y)\dy\,,
\]
and that moreover satisfies at least one of the following estimates,  for some $s\in (0,1]$:

\noindent - There exists a constant $c_{1,\Omega}>0$ such that for a.e. $x,y\in \Omega$\,:
\[\tag{K1}
0\le \K(x,y)\le c_{1,\Omega}\,|x-y|^{-(N-2s)}\,.
\]
- The second assumption concerns estimates involving the distance to the boundary. We assume that there exist constants $\gamma\in (0,1]$\,, $c_{0,\Omega},c_{1,\Omega}>0$ such that for a.e. $x,y\in \Omega$\,:
\[\tag{K2}
c_{0,\Omega}\,\p(x)\,\p(y) \le \K(x,y)\le \frac{c_{1,\Omega}}{|x-y|^{N-2s}}
\left(\frac{\p(x)}{|x-y|^\gamma}\wedge 1\right)
\left(\frac{\p(y)}{|x-y|^\gamma}\wedge 1\right)
\]
where we adopt the notation $\p(x):=\dist(x, \partial\Omega)^\gamma$.
Hypothesis (K2) introduces an exponent $\gamma$, which is a  characteristic of the operator and will play a big role in the results.

\medskip

\noindent\textbf{Remark. }(i) Defining an inverse operator $\AI$ implies that we are taking into account the Dirichlet boundary conditions. See more details in Section \ref{Sect.SFL+RFL+Examples}.

\noindent(ii) A key point of the hypothesis (K1) is the integrability property of $\K$. Indeed, by H\"older's inequality, the estimate implies that $\K(\cdot,x_0)\in \LL^r(\Omega)$ with $r<N/(N-2s)$\,, so that
\begin{equation}\label{K3.Holder}
\int_{\Omega} \psi(x)\,\K(x,x_0)\dx \le \|\psi\|_{r'}\|\K(\cdot,x_0)\|_r\,,\qquad\mbox{with }r'=\frac{r}{r-1}>\frac{N}{2s}\,.
\end{equation}
 
\noindent(iii) Hypotheses (K1) and (K2) are the only requirement on the operator $\A$\,, to obtain our main lower and upper estimates for the equation $u_t=\A \n(u)$\,, where $\n$  satisfies (N1). To be more precise, (A1),(A2) and (K1) imply the absolute upper bounds, and some smoothing effects. On the other hand, the stronger assumption (K2), implies weighted smoothing effects, and weighted $\LL^1$ estimates. We remark that in this paper we do not use the lower bound of Assumption (K2) to obtain sharp positivity estimates, a task that will be performed in \cite{BV-ppr2-2}, where we also derive sharp upper boundary behaviour and Harnack inequalities.

\noindent (iv)  Many formulas seen in the literature on these topics have convenient expressions in terms of the first eigenfunction.  If there exists a first eigenfunction of $\A$\,,  i.\,e., $\Phi_1\ge 0$ and $\A\Phi_1=\lambda_1\Phi_1$ for some $\lambda_1>0$, then hypothesis (K2) implies that $\Phi_1(x)\asymp \dist(x, \partial\Omega)^\gamma=\p(x)$\,, and we can rewrite $(K2)$ in the following equivalent form: \noindent There exist constants $\gamma\in (0,1]$\,, $c_{0,\Omega},c_{1,\Omega}>0$ such that for a.e. $x,y\in \Omega$\,:
\[\tag{K3}
c_{0,\Omega}\Phi_1(x)\Phi_1(y)\le \K(x,y)\le \frac{c_{1,\Omega}}{|x-x_0|^{N-2s}}
\left(\frac{\Phi_1(x)}{|x-y|^\gamma}\wedge 1\right)
\left(\frac{\Phi_1(y)}{|x-y|^\gamma}\wedge 1\right)\,.
\]
But we do not have to assume the existence of $\Phi_1$, even if our main examples do have a definite spectral sequence.
To cover the grater generality we use the weight $\p(x)=\dist(x, \partial\Omega)^\gamma$ in our estimates instead of $\Phi_1$\,, as it has been done in \cite{BV-PPR1}.

\noindent(v) The lower bound of assumption (K2) is weaker than the best known estimate on the Green function, both for the SFL and the RFL, which are the two main examples of fractional Laplacian operator on bounded domains. Indeed, for both mentioned operators we have
\[\tag{K4}
\K(x,y)\asymp \frac{1}{|x-x_0|^{N-2s}}
\left(\frac{\p(x)}{|x-y|^\gamma}\wedge 1\right)
\left(\frac{\p(y)}{|x-y|^\gamma}\wedge 1\right)\,.
\]
See \cite{Jak, Kul}, for a proof of the above estimates in the case of the RFL. For the SFL, the above bounds have been obtained separately: the bounds (K2) follow by the celebrated heat kernel estimates of Davies and Simon \cite{Davies1, DS} (for the case $s=1$), while the improved lower bound follows from the sharp lower heat kernel bounds of \cite{Zh2002}, see also \cite{BV-PPR1} or Section \ref{Sect.SFL+RFL+Examples} for more details.

We want to stress here that in order to obtain sharp boundary behaviour, we do not need (K4), but only the weaker (K2) estimates, that may hold for more general linear operators $\A$, as those mentioned in Section \ref{Sect.SFL+RFL+Examples}.

\medskip

\noindent {\bf Some notations. } The symbol $\infty$ will always denote $+\infty$. We also use the notation $a\asymp b$ if and only if there exist constants $c_0,c_1>0$ such that $c_0\,b\le a\le c_1 b$\,. We use the symbols $a\vee b=\max\{a,b\}$ and $a\wedge b=\min\{a,b\}$. \normalcolor We define the exponents for $\gamma\in [0,1]$:
\[
\vartheta_{i,\gamma}=\frac{1}{2s+(N+\gamma)(m_i-1)}\qquad\mbox{with }\qquad m_i=\frac{1}{1-\mu_i}>1
\]
that will appear in smoothing estimates.  Moreover, we denote by $\LL^p_{\p}(\Omega)$ the weighted $\LL^p$ space $\LL^p(\Omega\,,\, \p\dx)$, endowed with the norm
\[
\|f\|_{\LL^p_{\p}(\Omega)}=\left(\int_{\Omega} |f(x)|^p\p(x)\dx\right)^{\frac{1}{p}}\,.
\]
We will always consider bounded domains $\Omega$ with smooth boundary, at least $C^{2}$. The question of lower regularity of the boundary is not of concern here.

\section{The usual fractional Laplacian operators and other examples}\label{Sect.SFL+RFL+Examples}

We present now some of the nonlocal diffusion operators to which the previous scheme will be applied and precise conclusions derived. Firstly, we have been mainly interested in fractional Laplacian operators, hence we are concerned with the precise definition of such operators on a bounded domain of $\RR^N$. It is well-known that there is no ambiguity in the definition of the operators acting on the whole space, but this is not the case on bounded domains.

We have mentioned in the introduction the two main types of Fractional Laplacian operators with Dirichlet boundary conditions, namely RFL and SFL.
We first recall that both operators satisfy assumptions $(A1)$ and $(A2)$\,. Indeed in \cite{BSV2013} we provide a good functional setup both for the SFL and the RFL in the framework of fractional Sobolev spaces. Here we do not need such framework, but from the results that the reader can find there, it follows that assumptions $(A1), (A2)$ are satisfied. \normalcolor

\subsection{The Restricted Fractional Laplacian}\label{sect.RFL} On one hand, we can define a fractional Laplacian operator by using the integral representation of the whole space in terms of hypersingular kernels, namely
\begin{equation}\label{sLapl.Rd.Kernel}
(-\Delta_{\RR^N})^{s}  g(x)= c_{N,s}\mbox{
P.V.}\int_{\mathbb{R}^N} \frac{g(x)-g(z)}{|x-z|^{N+2s}}\,dz,
\end{equation}
where $c_{N,s}>0$ is a normalization constant, and ``restrict'' the operator to functions that are zero outside $\Omega$: we will denote the operator defined in such a way as $\A_2=(-\Delta_{|\Omega})^s$\,, and call it the \textit{restricted fractional Laplacian}\footnote{In the literature this is often called the fractional Laplacian on domains, but this simpler name may be confusing when the spectral fractional Laplacian is also considered, cf. \cite{BV-PPR1}}, RFL for short.
In this case, the initial and boundary conditions associated to the fractional diffusion equation \eqref{FPME.equation} read
\begin{equation}\label{FPME.Dirichlet.conditions.Restricted}
\left\{
\begin{array}{lll}
u(t,x)=0\,,\; &\mbox{in }(0,\infty)\times\RR^N\setminus \Omega\,,\\
u(0,\cdot)=u_0\,,\; &\mbox{in }\Omega\,,
\end{array}
\right.
\end{equation}
The boundary conditions can also be understood via the Caffarelli-Silvestre extension, see \cite{Caffarelli-Silvestre}, as explained in \cite{BSV2013}. The sharp form of the boundary behaviour for RFL has been investigated in \cite{RosSer}. We refer to \cite{BSV2013} for a careful construction of the RFL in the framework of fractional Sobolev spaces.  A probabilistic interpretation can be found for instance in \cite{BlGe}. In this case, assumptions (K1) and (K2) are satisfied  with exponent $\gamma=s$.  Indeed, we have that
\[
\K(x,y)\asymp \frac{1}{|x-y|^{N-2s}}
\left(\frac{\dist(x, \partial\Omega)}{|x-y|}\wedge 1\right)^s
\left(\frac{\dist(y, \partial\Omega)}{|x-y|}\wedge 1\right)^s
\]
and this has been proven in \cite{Jak, Kul}.

\subsection{Spectral Fractional Laplacian}\label{sect.SFL}

On the other hand, starting  from the classical Dirichlet Laplacian $\Delta_{\Omega}$ on the domain $\Omega$\,, then the so-called {\em spectral definition} of the fractional power of $\Delta_{\Omega}$ uses a formula in terms of the semigroup associated to the Laplacian, namely
\begin{equation}\label{sLapl.Omega.Spectral}
\displaystyle(-\Delta_{\Omega})^{s}
g(x)= \frac1{\Gamma(-s)}\int_0^\infty
\left(e^{t\Delta_{\Omega}}g(x)-g(x)\right)\frac{dt}{t^{1+s}}\,,
\end{equation}
which is equivalent to
\begin{equation}\label{sLapl.Omega.Spectral1}
\displaystyle(-\Delta_{\Omega})^{s}
g(x)=\sum_{j=1}^{\infty}\lambda_j^s\, \hat{g}_j\, \varphi_j(x)
\end{equation}
where  $\lambda_j>0$, $j=1,2,\ldots$ are the eigenvalues of the Dirichlet Laplacian on $\Omega$\,, written in increasing order and repeated according to their multiplicity, and $\varphi_j$ are the corresponding normalized eigenfunctions and
\[
\hat{g}_j=\int_\Omega g(x)\varphi_j(x)\dx\,,\qquad\mbox{with}\qquad \|\varphi_j\|_{\LL^2(\Omega)}=1\,.
\]
We will denote the operator defined in such a way as $\A_1=(-\Delta_{\Omega})^s$\,, and call it the \textit{spectral fractional Laplacian}, SFL for short. In this case, the initial and boundary conditions associated to the fractional diffusion equation \eqref{FPME.equation} read
\begin{equation}\label{FPME.Dirichlet.conditions.Spectral}
\left\{
\begin{array}{lll}
u(t,x)=0\,,\; &\mbox{in }(0,\infty)\times\partial\Omega\,,\\
u(0,\cdot)=u_0\,,\; &\mbox{in }\Omega\,.
\end{array}
\right.
\end{equation}
The boundary conditions can also be understood via the Caffarelli-Silvestre extension, see \cite{Cabre-Tan}, as explained in \cite{BSV2013}. For this operator, assumptions (K1) and (K2) are satisfied with $\gamma=s$. The estimates (K2) can be obtained by the following Heat kernel estimates $H(t,x,y)$,  for the case $s=1$
\begin{equation}\label{Heat.Kernel.est}
H(t,x,x_0)\asymp \left(\frac{\varphi_1(x)}{|x-x_0|}\wedge 1\right)
\left(\frac{\varphi_1(x_0)}{|x-x_0|}\wedge 1\right)\,\frac{\ee^{-\frac{c_2|x-x_0|^2}{t}}}{t^{d/2}}\,.
\end{equation}
with the help of the formula
\begin{equation}\label{green.heat}
\K(x,y)=\int_0^{\infty}\frac{H(t,x,y)}{t^{1-s}}\dt\,.
\end{equation}
The upper bounds for the Heat kernel \eqref{Heat.Kernel.est} can be found in \cite{D1, D2, Davies1, Davies2, DS}\,, while the lower bounds have been obtained later in \cite{Zh2002}\,.


\subsection{More examples}\label{sect.MoreExamples}
Our general framework applies to a number of important operators. Here below we list  some that we find relevant.

\noindent\textbf{Censored fractional Laplacian and operators with general kernels. }Censored stochastic processes have been introduced by Bogdan et al. in \cite{bogdan-censor}. Very roughly speaking, they correspond to processes with a different homogeneous Dirichlet condition on the boundary of the domain $\Omega$; it describes the situation when a particle that hits the boundary is forbidden to jump outside but is allowed to continue its path if it happens to continue inside $\Omega$.  The infinitesimal operator has the form \cite{bogdan-censor, Song-coeff}
\[
\A f(x)=\mathrm{P.V.}\int_{\Omega}\left(f(x)-f(y)\right)\frac{a(x,y)}{|x-y|^{N+2s}}\dy\,,\qquad\mbox{with }\frac{1}{2}<s<1\,,
\]
where $a(x,y)$ is a measurable symmetric function bounded between two positive constants, satisfying some further assumptions; a sufficient assumption is $a\in C^{1}(\overline{\Omega}\times\overline{\Omega})$. Whenever $s\in \left(\frac{1}{2}, 1\right)$, the Green function $\K(x,y)$ of $\A$ satisfies the strongest assumption $(K_4)$ with $\gamma=s-1/2$\,, namely
\begin{equation*}
\K(x,y)\asymp \frac{1}{|x-y|^{N-2s}}
\left(\frac{\p(x)}{|x-y|^\gamma}\wedge 1\right)
\left(\frac{\p(y)}{|x-y|^\gamma}\wedge 1\right)\,, \qquad\mbox{with }\gamma=s-\frac{1}{2}\,.
\end{equation*}
This bounds has been proven in Corollary 1.2 of  \cite{Song-coeff}, as a consequence of sharp heat kernel estimates.

\medskip

\noindent\textbf{Fractional operators with general kernels. }Consider integral operators of the following form\vspace{-1mm}
\[
\A f(x)=\mathrm{P.V.}\int_{\RR^N}\left(f(x+y)-f(y)\right)J(x,y)\dy\,.
\]
where $J(x,y)$ is a kernel such that\vspace{-1mm}
\[
J(x,y)=\frac{K(x,y)}{|x-y|^{N+2s}}
\]
and $K$ is a measurable symmetric function bounded between two positive constants, satisfying
\[
\big|K(x,y)-K(x,x)\big|\,\chi_{|x-y|<1}\le c |x-y|^\sigma\,,\qquad\mbox{with }0<s<\sigma\le 1\,,
\]
for some positive $c>0$. We can allow even more general kernels. Take for instance a function $\eta: [0,+\infty)\to [0,+\infty)$
\[
\eta(r)= 1\quad \mbox{if }0\le r\le 1\qquad\mbox{and}\qquad c_1\ee^{c_2\,r^\beta}\le \eta(r)\le c_3\ee^{c_4\,r^\beta}\quad\mbox{if }r>1\,,
\]
for some $c_1,\dots,c_4>0$ and $\beta\in [0,\infty]$\,, (the case above was for $\beta=0$) and construct the kernel\vspace{-1mm}
\[
J(x,y)=K(x,y) \overline{J}(|x-y|)
=\frac{K(x,y)}{|x-y|^{N+2s}}\left\{
\begin{array}{ll}
\eta(|x-y|)^{-1}& ~~\mbox{if }0\le \beta<\infty\,,\\
\chi_{|x-y|<1} & ~~\mbox{if }\beta=\infty\,.\\
\end{array}
\right.
\]
where $K$ and $\eta$ are as above. For all $s\in (0, 1]$, the Green function $\K(x,y)$ of $\A$ satisfies the strongest assumption $(K_4)$ with $\gamma=s$\,, namely\vspace{-1mm}
\begin{equation}\label{K4-examples}
\K(x,y)\asymp \frac{1}{|x-y|^{N-2s}}
\left(\frac{\p(x)}{|x-y|^\gamma}\wedge 1\right)
\left(\frac{\p(y)}{|x-y|^\gamma}\wedge 1\right)\,, \qquad\mbox{with }\gamma=s\,.
\end{equation}
This bounds has been proven in Corollary 1.4 of \cite{Kim-Coeff}, as a consequence of sharp heat kernel estimates.

\medskip

\noindent\textbf{Spectral powers of uniformly elliptic operators. }Consider a linear operator $A$ in divergence form:
\[
A=\sum_{i,j=1}^N\partial_i(a_{ij}\partial_j)\,,
\]
with bounded measurable coefficients, which are uniformly elliptic. The uniform ellipticity allows to build a self-adjoint operator on $\LL^2(\Omega)$ with discrete spectrum $(\lambda_k, \phi_k)$\,. Using the spectral theorem, we can construct the spectral power of such operator, defined as follows:
\[
\A f(x):=A^s\,f(x):=\sum_{k=1}^\infty \lambda_k^s \hat{f}_k \phi_k(x)\qquad\mbox{where }\qquad \hat{f}_k=\int_\Omega f(x)\phi_k(x)\dx\,.
\]
We refer to the books \cite{Davies1,Davies2} for further details. This construction is similar to the SFL, formally replacing the Laplacian with a more general uniformly elliptic operator. The Green function satisfies the estimates $(K2)$ and $(K3)$, with $\gamma=1$\,, as a consequence of the heat kernel bounds of the form \eqref{Heat.Kernel.est} and formula \eqref{green.heat}\,.  Indeed, since there exist a first eigenfunction $\Phi_1(x)\asymp \dist(x, \partial\Omega)$\,, we can rewrite $(K2)$ in the equivalent form (K3):
\[
c_{0,\Omega}\Phi_1(x)\Phi_1(y)\le \K(x,y)\le \frac{c_{1,\Omega}}{|x-x_0|^{N-2s}}
\left(\frac{\Phi_1(x)}{|x-y|}\wedge 1\right)
\left(\frac{\Phi_1(y)}{|x-y|}\wedge 1\right)\,.
\]
As a general statement, we can treat the class of intrinsically ultra-contractive operators introduced in \cite{DS}, cf. also \cite{Song-intrinsic-ultra, Davies1,Davies2}. In \cite{Song-subordinate1,Song-subordinate2} analogous estimates are derived for nontrivial perturbation of the above mentioned operators.

On the other hand, it is possible to consider more general classes of linear operators, for instance with singular or degenerate coefficients, and their powers, see \cite{Davies1,Davies2,DS}, but in this case the power $\gamma\in (0,1]$ may change depending on the singularity/degeneracy of the coefficients. Our theory does not cover such cases.

\medskip

\noindent\textbf{Sums of two fractional operators. } Operators of the form
\[
\A=(\Delta_{|\Omega})^{s}+(\Delta_{|\Omega})^{\sigma}\,,\qquad\mbox{with }0<\sigma<s\le 1\,,
\]
where $(\Delta_{|\Omega})^{s}$ is the RFL defined in section \eqref{sect.RFL}\,. The Green function $\K(x,y)$ of $\A$ satisfies the strongest assumption $(K4)$ in the form \eqref{K4-examples}, with $\gamma=s$\,. The above bounds are consequence of sharp heat kernel estimates, which has been proven in \cite{Song-Sum3} Corollary 1.2, for the case $0<\sigma<s<1$\,. The limit case $s=1$ and $\sigma\in (0,1)$ has been proven in \cite{Song-Sum2} Theorem 1.1\,, in which case we have $\gamma=s=1$.
 
\medskip

\noindent\textbf{Sum of the Laplacian and operators with general kernels. }In the case
\[
\A=a\Delta+ A_s \,,\qquad\mbox{with }0<s< 1\quad\mbox{and}\quad a\ge 0\,,
\]
where
\[
A_sf(x)=\mathrm{P.V.}\int_{\RR^N}\left(f(x+y)-f(y)-\nabla f(x)\cdot y \chi_{|y|\le 1}\right)\chi_{|y|\le 1}\rd\nu(y)\,.
\]
where the measure $\nu$ on $\RR^N\setminus\{0\}$ is invariant under rotations around origin and satisfies
\[
\int_{\RR^N} (1\vee|x|^2)\,\rd\nu(y)<\infty\,.
\]
More precisely $\rd\nu(y)= j(y)\dy$ with $j:(0,\infty)\to(0,\infty)$ is given by
\[
j(r):=\int_0^\infty\frac{\ee^{r^2/(4t)}}{(4\pi\,t)^{N/2}}\rd\mu(r)\qquad\mbox{with}\qquad \int_0^\infty (1\vee t)\rd\mu(t)<\infty\,.
\]
In Theorem 1.4, of \cite{Song-Sum1Gen} the authors prove that the Green function $\K(x,y)$ of $\A$ satisfies the strongest assumption $(K_4)$ in the form \eqref{K4-examples}, with $s=1$ and $\gamma=1$\,.
 
\medskip

\noindent\textbf{Schr\"odinger equations for non-symmetric diffusions. }In the case
\[
\A=A + \mu\cdot\nabla + \nu\,,
\]
where $A$ is uniformly elliptic and is allowed to take both divergence and non-divergence form:
\[
A_1=\frac{1}{2}\sum_{i,j=1}^N\partial_i(a_{ij}\partial_j) \qquad\mbox{or}\qquad A_2=\frac{1}{2}\sum_{i,j=1}^N a_{ij}\partial_{ij}\,,
\]
and we assume $C^1$ coefficient $a_{ij}$, uniformly elliptic. Finally, $\mu,\nu$ are measures belonging to suitable Kato classes, we refer to \cite{song-NonSymm-1} for more details. \normalcolor The Green function $\K(x,y)$ of $\A$ satisfies the strongest assumption $(K_4)$ in the form \eqref{K4-examples}, with $s=1$ and $\gamma=1$\,, cf. \cite{song-NonSymm-1}.
 
\medskip

\noindent\textbf{Gradient perturbation of restricted fractional Laplacians. }In the case
\[
\A=(\Delta_{|\Omega})^{s}+b\cdot\nabla
\]
where $b$ is a vector valued function belonging to a suitable Kato class, we refer to \cite{Song-Drift} for further details. The Green function $\K(x,y)$ of $\A$ satisfies the strongest assumption $(K_4)$ in the form \eqref{K4-examples}, with $\gamma=s$\,, see  Corollary 1.4 of \cite{Song-Drift}.
 
\medskip

\noindent\textbf{Relativistic stable processes. }In this case
\[
\A=c-\left(c^{1/s}-\Delta\right)^s\,,\qquad\mbox{with }c>0\,,\mbox{ and }0<s\le 1\,.
\]
For more details about the associated relativistic stable process we refer to \cite{Song-Rel}. The Green function $\K(x,y)$ of $\A$ satisfies the strongest assumption $(K4)$ in the form \eqref{K4-examples}, with $\gamma=s$, see Theorem 1.3 of \cite{Song-Rel}.

\medskip

In the selection of these examples we are indebted to a number of authors who have provided the basic estimates, and we would like to specially mention Prof. Renming Song.

%
%
%
\section{Weak dual solutions, existence and uniqueness}

We are going to find solutions of Equation \eqref{FPME.equation} with zero Dirichlet boundary data in the following generalized sense called weak dual solution, that  is expressed in terms of the a problem involving the inverse of the operator $\A$\,. We recall that the inverse $\AI$ contains the information on the boundary data.
\begin{defn}\label{Def.Very.Weak.Sol.Dual} A function $u$ is a {\sl weak dual} solution to the Dirichlet Problem for Equation \eqref{FPME.equation} in $Q_T=(0,T)\times \Omega$ if:
\begin{itemize}
\item $u\in C((0,T): \LL^1_{\p}(\Omega))$\,, $\n(u) \in \LL^1\left((0,T):\LL^1_{\p}(\Omega)\right)$;
\item  The identity
\begin{equation}
\displaystyle \int_0^T\int_{\Omega}\AI (u) \,\dfrac{\partial \psi}{\partial t}\,\dx\dt
-\int_0^T\int_{\Omega} \n(u)\,\psi\,\dx \dt=0.
\end{equation}
holds for every test function $\psi$ such that  $\psi/\p\in C^1_c((0,T): \LL^\infty(\Omega))$\,.
\end{itemize}\end{defn}
More precisely, we will solve the problem consisting of Equation \eqref{FPME.equation} with homogeneous Dirichlet conditions plus given initial data. We will call this problem (CDP).

\begin{defn}\label{Def.Very.Weak.Sol.Dual.CDP}
A {\sl weak dual} solution to the Cauchy-Dirichlet problem \textbf{(CDP)} is a  weak dual solution to Equation \eqref{FPME.equation} such that moreover $u\in C([0,T): \LL^1_{\p}(\Omega))$  and $u(0,x)=u_0\in \LL^1_{\p}(\Omega)$.
\end{defn}
This kind of solution has been first introduced in \cite{BV-PPR1}.
Roughly speaking, we are considering the weak solution to the ``dual equation'' $\partial_t U=- u^m$\,, where $U=\AI u$\,, posed on the bounded domain $\Omega$ with homogeneous Dirichlet conditions. The dual equation has been used in the case $s=1$ by Pierre \cite{Pierre}, to prove uniqueness of solutions with measure initial data. In the fractional setting, it is used in \cite{Vaz2012} for the RFL on the whole $\RR^N$.

\noindent Two useful remarks. (i) The finite existence time $T>0$ is used in the definition for generality, in view of other possible applications. Below, the maximal choice $T=\infty$ is always used.

\noindent (ii) Notice that the condition on $\psi$ forces it to vanish on the boundary. Actually,  $\psi/\p\in C^1_c((0,T): \LL^\infty(\Omega))$ implies that $\|\psi(t,\cdot)/\p\|_{\LL^\infty(\Omega)}<+\infty$  and $\|\partial_t\psi(t,\cdot)/\p\|_{\LL^\infty(\Omega)}<+\infty$ for all $t\ge 0$ and moreover, are compactly supported functions of $t>0$ therefore in $\LL^1(0,\infty)$.\normalcolor

In the course of our study we will need a somewhat more restricted class of solutions.
\begin{defn}\label{Def.Very.Weak.Sol.Dual.2} We consider a class $\mathcal{S}_p$ of nonnegative  weak dual solutions $u$ to (CDP) with initial datum in $u_0\in\LL^1_{\p}(\Omega)$\,, such that  (i) the map $u_0\mapsto u(t)$ is  order preserving in $ \LL^1_{\p}(\Omega)$; (ii) for all $t>0$ we have $u(t)\in L^p(\Omega)$ for some $p\ge 1$.
\end{defn}
This class is not empty, and indeed we will show in Corollary \ref{Prop.VWSp-WDS} that it contains the semigroup (mild) solutions with initial data in $\LL^p(\Omega)$\,. The strategy that we will use is the following: first we show in Proposition \ref{Prop.Mild-WDS} that semigroup (mild) solutions are weak dual solutions. Next we show, in Proposition \ref{Prop.Lp.Mild} that semigroup solutions are stable in $\LL^p$. The combination of the two mentioned propositions implies that mild solutions with $u_0\in\LL^p(\Omega)$ are weak dual weak solutions belonging to the class $\mathcal{S}_p$\,. This will be done in Section \ref{Sect.Various.Sols}.   We advance the results for the reader's convenience.

\begin{thm}[Existence of weak dual solutions]\label{thm.L1weight.exist}
For every  nonnegative $ u_0\in\LL^1_{\p}(\Omega)$ there exists a minimal weak dual solution to the $(CDP)$. Such a solution is obtained as the monotone limit of the semigroup (mild) solutions that exist and are unique. The minimal  weak dual solution is continuous in the weighted space $u\in C([0,\infty):\LL^1_{\p}(\Omega))$. Mild solutions are weak dual solutions and the set of such solutions has the properties needed to form a class of type~$\mathcal{S}$.
\end{thm}
\begin{thm}[Uniqueness]\label{thm.L1weight.uniq}
The solution constructed in Theorem $\ref{thm.L1weight.exist}$ by approximation of the initial data from below is unique. We call it the minimal solution. In this class of solutions the standard comparison result holds, and also the estimates of Proposition $\ref{thm.L1weight.contr.psi}$\,.
\end{thm}
We postpone the proof of these results because it relies on the a priori estimates for weak dual solutions that we will derive next and are the quantitative basis of the paper.
The proof of existence will be then done in Section \ref{sect.exist.uniq.WDS}, since it depends on the weighted $\LL^1$ estimates of Section \ref{Sect.Weighted.L1}, which in turn depend on the upper estimates that we state in the next section.
%
%
\section{The two first estimates}\label{sect.frist.two}
\subsection{Pointwise estimates for weak dual solutions}\label{sect.dualweak.estimates}
We begin by proving the basic pointwise estimates, needed to prove both the absolute upper bounds and the smoothing effects.
\begin{prop}\label{prop.point.est}  let $u\ge 0$ be a solution in the class $\mathcal{S}_p$ of very weak solutions to Problem $(CDP)$  with $p>N/2s$. Then,
\begin{equation}\label{thm.NLE.PME.estim.0}
\int_{\Omega}u(t,x)\K(x , x_0)\dx\le \int_{\Omega}u_0(x)\K(x , x_0)\dx \qquad\mbox{for all $t> 0$\,.}
\end{equation}
Moreover, for almost every $0< t_0\le t_1 \le t$ and almost every $x_0\in \Omega$\,, we have
\begin{equation}\label{thm.NLE.PME.estim}
\left(\frac{t_0}{t_1}\right)^{\frac{1}{\mu_0}}(t_1-t_0)\,\n(u(t_0,x_0))
\le \int_{\Omega}\big[u(t_0,x)-u({t_1},x)\big]\K(x , x_0)\dx  \le (m_0-1)\frac{t^{\frac{1}{\mu_0}}}{t_0^{\frac{1-\mu_0}{\mu_0}}} \,\n(u(t,x_0))\,.
\end{equation}
\end{prop}
\noindent\textsc{Sketch of the proof. } We adapt the scheme of the proof of Proposition 4.2 of \cite{BV-PPR1}, in which we have treated the case $\n(u)=u^m$.  We provide a formal proof here  in order to understand the main ideas. The fully rigorous proof is given in Appendix \ref{app.proof.5.1}.

The first issue is to take a test function of the form $\psi(t,x)=\psi_1(t)\psi_2(x)$, with $\psi_1(\tau)=\chi_{[t_0,t_1]}(\tau)$ and $\psi_2(x)=\K(x_0,x)$\,, i.e. to be able to admit the Green function as a test. Of course, this choice of test function is not directly admissible in Definition \ref{Def.Very.Weak.Sol.Dual} of weak dual solution, therefore it has to be justified though a careful approximation, see the proof in Appendix \ref{app.proof.5.1}. Plugging such test function in the definition of weak dual solution gives the identity
\begin{equation}\label{step.1.thm.repr.00}\begin{split}
\int_{\Omega}u(t_0,x) \K(x_0,x)\dx - \int_{\Omega}u({t_1},x)\K(x_0,x)\dx
=\int_{t_0}^{t_1} \n(u(\tau,x_0))\rd\tau\,.
\end{split}
\end{equation}
Next, we use the monotonicity estimates \eqref{CP.Monotonicity} of \cite{CP-JFA}, namely that the function
$ t\mapsto t^{\frac{1}{\mu_0}}\,\n(u(t,x))$  is nondecreasing in $t>0$ for a.e. $x\in \Omega$\,, recalling that $\frac{1}{\mu_0}=\frac{m_0}{m_0-1}$, to get for all $0\le t_0\le t_1 \le t$:
\begin{equation}\label{step.1.thm.repr.01}
\begin{split}
\left(\frac{t_0}{t_1}\right)^{\frac{m_0}{m_0-1}}(t_1-t_0)\n(u(t_0,x_0))
&\le \int_{t_0}^{t_1}\n(u(\tau,x_0)) \rd\tau \le \frac{m_0-1}{t_0^{\frac{1}{m_0-1}}} t^{\frac{m_0}{m_0-1}}\,\n(u(t,x_0))\,.
\end{split}
\end{equation}
Joining \eqref{step.1.thm.repr.00} and \eqref{step.1.thm.repr.01} gives \eqref{thm.NLE.PME.estim}\,.\qed

\noindent\textbf{Remarks. }(i) A solution  $u(x,t)\ge 0$ belonging to the class $\mathcal{S}_p$ with $p>N/2s$ has the property that $u(t)\in \LL^p(\Omega)$, therefore $\int_{\Omega}u(t,x)\K(x , x_0)\dx<+\infty$ for all $t>0$\,, simply as a consequence of H\"older inequality, as in \eqref{K3.Holder}. Therefore, using the pointwise inequality \eqref{thm.NLE.PME.estim} we may conclude that $u(\cdot,t)$ is bounded for all $t>0$. On the other hand, we only assume that  $0\le u_0\in \LL^1_{\p(\Omega)}$, hence $\int_{\Omega}u_0(x)\K(x , x_0)\dx$ may be unbounded.

\noindent (ii) The most relevant part of the last estimate for the applications  is represented by the middle term: we will use it to deduce both lower and upper estimates, that will be combined into quantitative Harnack inequalities; this will be done in a forthcoming paper \cite{BV-ppr2-2}, while the case $\n(u)=|u|^{m-1}u$ has been treated in \cite{BV-PPR1}.


\subsection{Absolute upper bounds }

We use the pointwise lower  estimates of Proposition \ref{prop.point.est} as a  basis to prove a number of quantitative estimates. The first one is the absolute upper bound that is formulated in terms of the nonlinearity $\n$ and its Legendre transform $\nl$\,, and has an interesting intrinsic form. \normalcolor

\begin{thm}[Absolute upper estimate]\label{thm.Upper.PME}
Let  $u$ be a nonnegative weak dual solution corresponding to $u_0\in \LL^1_{\p}(\Omega)$. Then, there exists universal constants $K_0, K_1, K_2>0$ such that the following estimates hold true for all $t> 0$\,:
\begin{equation}\label{thm.Upper.PME.Absolute.F}
\n\left(\|u(t)\|_{\LL^\infty(\Omega)}\right)
\le \nl\left(\frac{K_1}{t}\right)\,.
\end{equation}
Moreover, there exists a time $\tau_1(u_0)$ with $0\le \tau_1(u_0)\le K_0$ such that $\|u(t)\|_{\LL^\infty(\Omega)}\le 1$ for all $t\ge \tau_1$ and
\begin{equation}\label{thm.Upper.PME.Absolute}
\|u(t)\|_{\LL^\infty(\Omega)}\le\frac{K_2}{t^{\frac{1}{m_i-1}}}\qquad\mbox{with $i=0$ if $t\le K_0$ and $i=1$ if $t\ge K_0$}\,.
\end{equation}
\end{thm}
The Legendre transform of $\n$ is defined as a function $\nl:\RR\to\RR$ with
\begin{equation}\label{def.F*}
\nl(z)=\sup_{r\in\RR}\big(zr-\n(r)\big)=z\,(\n')^{-1}(z)-\n\left((\n')^{-1}(z)\right)
=\n'(r)\,r + \n(r)\,,
\end{equation}
with the choice $r=(\n')^{-1}(z)$\,.

\noindent\textbf{Remarks. }(i) This absolute bound proves a strong regularization which is independent of the initial datum. While the first inequality is intrinsic and accurate, the second is just convenient, so as to have more explicit estimate. The constants $K_0,K_1, K_2>0$ depend only on $N, s, \mu_0,\mu_1$ and $\Omega$\,, but not on $u$\,, and have an explicit form given in the proof.

\noindent(ii) Let us briefly explain  the role of $K_0$: roughly speaking, it is an upper estimate of the time $\tau_1(u_0)$ after which $\|u(t)\|_{\LL^\infty(\Omega)}\le 1$, which is important because the behaviour of $u\mapsto \n(u)$ may change when $u\le 1$ or when $u\ge 1$: the powers in the above estimates change in the two cases in a quite sharp way.

\noindent(iii) The bounds \eqref{thm.Upper.PME.Absolute} hold for any $t> 0$\,, but they are not sharp for small times; precise bounds for small times are the smoothing effects of Theorem \ref{thm.Upper.1.PME}, and they depend on the initial datum.

\noindent(iv) This result will also imply the sharp upper boundary behaviour, see the forthcoming paper \cite{BV-ppr2-2}.

Before proceeding wit the proof of this result, we recall first some basic facts about the Legendre transform. More details can be found in Appendix \ref{App.F.N1}\,.

\subsection{Properties of $\n$ and its Legendre transform}\label{ssec.Lecgendre}
Let us recall that hypothesis (N1) on $\n$ implies its convexity, and it is equivalent to (N2) (see Appendix \ref{App.F.N1}) which we recall here
\begin{equation}\label{hyp.1.nonlinearity.N1+N3-recall}
\mu_0\le \frac{\n(r)\n''(r)}{[\n'(r)]^2}\le \mu_1 \qquad\mbox{or}\qquad
1-\mu_1\le \left(\frac{\n(r)}{\n'(r)}\right)'\le 1-\mu_0
\qquad\mbox{a.e. }r>0\,.
\end{equation}
We can integrate the latter to get (recall that $0<\mu_0\le \mu_1<1$):
\begin{equation}\label{hyp.N3.integrated.0}
(1-\mu_1)r\le \frac{\n(r)}{\n'(r)}\le (1-\mu_0)r\qquad\mbox{or}\qquad
\frac{1}{(1-\mu_0)r}\le \frac{\n'(r)}{\n(r)}\le \frac{1}{(1-\mu_1)r}
\end{equation}
which can be rewritten in the form
\begin{equation}\label{hyp.N3.integrated}
\mu_0 r\,\n'(r)\le r\,\n'(r)-\n(r)\le \mu_1 r\,\n'(r)
\end{equation}
By definition \eqref{def.F*} of Legendre transform, it is easy to show that
\begin{equation}\label{def.F*.deriv}
\nl'(z)=(\n')^{-1}(z)=r\qquad\mbox{and}\qquad \nl''(z)=\frac{1}{\n''\left((\n')^{-1}(z)\right)}=\frac{1}{\n''(r)}\ge 0
\end{equation}
As a consequence of the above inequalities we obtain:
\[
\mu_0 z \nl'(z)=\mu_0 r\,\n'(r)\le \nl(z)=\n'(r)\,r + \n(r)\le \mu_1 r\,\n'(r)= \mu_1 z \nl'(z)
\]
from which we derive
\begin{equation}\label{def.F*.2}
\mu_0\,z\le\frac{\nl(z)}{\nl'(z)}\le \mu_1\,z\qquad\mbox{or}\qquad
\mu_0\le\frac{1}{z}\left(\frac{\nl(z)}{\nl'(z)}-\frac{\nl(0)}{\nl'(0)}\right)\le \mu_1
\end{equation}
since $\nl(0)/\nl'(0)=0$. This finally gives the ``dual'' (N1) hypotheses (in an equivalent form):
\begin{equation}\label{def.F*.3}
\mu_0\le\left(\frac{\nl(z)}{\nl'(z)}\right)'\le \mu_1\qquad\mbox{or}\qquad 1-\mu_1\le \frac{\nl(z)\nl''(z)}{[\nl'(z)]^2} \le 1-\mu_0
\end{equation}

Finally, we recall a scaling property of the Legendre transform: for any $\varepsilon>0$, letting $\n_\varepsilon(r)=\varepsilon\n(r)$ we get $\nl_\varepsilon(z)=\varepsilon\nl(z/\varepsilon)$\,. This scaling property, together with the Young inequality, gives
\begin{equation}\label{Young.F}
a\,b\, \le \n_\varepsilon(a)+\nl_\varepsilon(b)=\varepsilon \n(a) +\varepsilon \nl\left(\frac{b}{\varepsilon}\right)\qquad\mbox{for all $a,b\ge 0$\,.}
\end{equation}

\medskip

Thanks to hypothesis (N1) it is possible to understand the behaviour at zero and at infinity of $\n$ and $\nl$, see Lemma \ref{Lem.N1.F} in Appendix \ref{App.F.N1} for a proof. We summarize the relevant information about $\n$ and $\nl$ in the table below.   We recall that $m_i=1/(1-\mu_i)$ or, equivalently, $\mu_i=(m_i-1)/m_i$\,.
\begin{center}
\begin{tabular}{| c | c | c |}
\hline & & \\
    $\n$
        & /
        & $\nl$ \\
& & \\ \hline & & \\
$ 1-\mu_1\le \left(\dfrac{\n(r)}{\n'(r)}\right)'\le 1-\mu_0$
        & $(N1)$
        & $\mu_0 \le\left(\dfrac{\nl(z)}{\nl'(z)}\right)'\le \mu_1$ \\
& & \\ \hline & & \\
        $\mu_0\le \dfrac{\n(r)\n''(r)}{[\n'(r)]^2}\le \mu_1$
        & $(N2)$
        & $1-\mu_1\le \dfrac{\nl(z)\nl''(z)}{[\nl'(z)]^2} \le 1-\mu_0$ \\
& & \\ \hline & & \\
    $\kb \left(\dfrac{r}{r_0}\right)^{m_1}\le \dfrac{\n(r)}{\n(r_0)} \le \ka \left(\dfrac{r}{r_0}\right)^{m_0}$
        & for all $0\le r\le r_0$
        & $\kb \left(\dfrac{r}{r_0}\right)^{\frac{m_0}{m_0-1}}\le \dfrac{\nl(r)}{\nl(r_0)}\le \ka \left(\dfrac{r}{r_0}\right)^{\frac{m_1}{m_1-1}}$ \\
& & \\ \hline & & \\
   $ \left(\dfrac{r}{r_0}\right)^{m_0}\le \dfrac{\n(r)}{\n(r_0)} \le \left(\dfrac{r}{r_0}\right)^{m_1}$
        & for all $r\ge r_0\ge 0$
        & $ \left(\dfrac{r}{r_0}\right)^{\frac{m_1}{m_1-1}}\le \dfrac{\nl(r)}{\nl(r_0)}\le\left(\dfrac{r}{r_0}\right)^{\frac{m_0}{m_0-1}}$ \\
& & \\ \hline
\end{tabular}
\end{center}


\subsection{Proof of the absolute upper bounds (AUB)}

The proof is a consequence of estimates (K1) for the Green function. \normalcolor On the other hand, the strategy of the proof follows  the methods of \cite{BV-PPR1}, in which we have treated the case $\n(u)=u^m$ with $m>1$. We first need a lemma.\normalcolor

\begin{lem}[Integral Green function estimates]\label{Lem.Green}Let $\K$ be the Green function of $\A$. Then, the $(K1)$ estimates imply that there exist a constant $c_{2,\Omega}(q)>0$ such that
\begin{equation}\label{Lem.Green.est.Upper.I}
\sup_{x_0\in\Omega}\int_{\Omega}\K^q(x , x_0)\dx \le c_{2,\Omega}(q)\qquad\mbox{for all $0<q<\dfrac{d}{d-2s}$\,,}
\end{equation}
\end{lem}

\noindent {\sl Proof.~}The proof follows by H\"older's inequality. A stronger version of this lemma will be proved in \cite{BV-ppr2-2}.\qed

\noindent We now proceed with the proof of Theorem \ref{thm.Upper.PME}.

\noindent$\bullet~$\textsc{Step 0. }\textit{Reduction. }We begin by observing that it is not restrictive to assume that $u\in \mathcal{S}_p$ with $p>N/(2s)$. Indeed, consider a weak dual solution corresponding to the initial datum $0\le u_{0,n}\in \LL^1_{\p}(\Omega)$. Define the sequence $0\le u_{0,n}:=n\wedge u_0\in \LL^\infty(\Omega)$. We know that $u_{0,n}$ converges monotonically to $u_{0,n}\in \LL^1_{\p}(\Omega)$ in the strong $\LL^1_{\p}$ topology. Let $u_n$ be the unique mild solutions constructed in Theorem \ref{thm.CP-JFA}, corresponding  to $u_{0,n}$ as initial data, which we know to belong to the class $\mathcal{S}_p$ thanks to the results of Section \eqref{Sect.Various.Sols}. For such solutions inequality \eqref{thm.Upper.PME.Absolute.F} holds true, hence
\begin{equation*}
\n\left(\|u_n(t)\|_{\LL^\infty(\Omega)}\right)
\le \nl\left(\frac{K_1}{t}\right)\,.
\end{equation*}
Taking the limit in the above inequality, using the lower semicontinuity of the $\LL^\infty$ norm, gives the desired result for the limit solution $u$\,, which by uniqueness - see Theorem \ref{thm.L1weight.uniq} - turns out to be the minimal weak dual solution we started with.

\noindent$\bullet~$\textsc{Step 1. }\textit{Fundamental upper estimates. }We first recall the lower pointwise estimate of Proposition \ref{prop.point.est}, that holds for any solution $u\in \mathcal{S}_p$ with $p>N/(2s)$:  for all $0\le t_0\le t_1 $ and $x_0\in \Omega$\,, we have that
\begin{equation}\label{Upper.PME.Step.1.1}
\left(\frac{t_0}{t_1}\right)^{\frac{1}{\mu_0}}(t_1-t_0)\,\n(u(t_0,x_0)) \le \int_{\Omega}u(t_0,x)\K(x , x_0)\dx - \int_{\Omega}u({t_1},x)\K(x , x_0)\dx\,.
\end{equation}
We choose $t_1=2t_0$ and recall that $u\ge0$\,, so that the above inequality \eqref{Upper.PME.Step.1.1} implies that
\begin{equation}\label{Upper.PME.Step.1.2}
\n(u(t_0,x_0)) \le \frac{2^{\frac{1}{\mu_0}}}{t_0}\int_{\Omega}u(t_0,x)\K(x , x_0)\dx\qquad\mbox{for all $t_0> 0$ and $x_0\in \Omega$}\,.
\end{equation}
This is a fundamental upper bound which encodes both the smoothing effect and the absolute upper bound, and therefore it is sharp both for large and for small times. A remarkable aspect of this upper bound is that it compares the $\LL^\infty$ norm and some integral norms, both at the same time $t_0> 0$\,.

\noindent The fact that $u\in \mathcal{S}_p$ guarantees that  $u(t,\cdot)\in \LL^p(\Omega)$ for all $t> 0$\,, with $p>N/(2s)$\,, so that
\[
\int_{\Omega}u(t_0,x)\K(x , x_0)\dx \le \|u(t_0)\|_{\LL^p(\Omega)}\,\|\K(\cdot , x_0)\|_{\LL^q(\Omega)}\le c_{2,\Omega}(q)\|u(t_0)\|_{\LL^p(\Omega)}<+\infty
\]
because $\K(\cdot , x_0)\in\LL^q(\Omega) $ for all $0< q<N/(N-2s)$\,, see Lemma \ref{Lem.Green}, in particular for $q=(p-1)/p$. Therefore, we have
\begin{equation}\label{Upper.PME.Step.1.3}
\n(u(t_0,x_0)) \le c_{2,\Omega}(q) \frac{2^{\frac{1}{\mu_0}}}{t_0} \|u(t_0)\|_{\LL^p(\Omega)}\qquad\mbox{for all $t_0> 0$ and $x_0\in \Omega$}\,.
\end{equation}
so that $u(t_0)\in \LL^\infty(\Omega)$ for all $t_0>0$.

As a byproduct of inequality \eqref{Upper.PME.Step.1.3}, we have proved that the class $\mathcal{S}_p$ with $p>N/(2s)$ and the class $\mathcal{S}_\infty$ are the same.

\noindent$\bullet~$\textsc{Step 2. }\textit{Proof of inequality $\eqref{thm.Upper.PME.Absolute.F}$. }Let us estimate the right-hand side of the fundamental upper bound \eqref{Upper.PME.Step.1.2} in another way as follows:
\begin{equation}\label{Upper.PME.Step.2.0}
\n(u(t_0,x_0)) \le \frac{2^{\frac{1}{\mu_0}}}{t_0}\int_{\Omega}u(t_0,x)\K(x , x_0)\dx
\le \|u(t_0)\|_{\LL^\infty(\Omega)}\frac{2^{\frac{1}{\mu_0}}}{t_0}\int_{\Omega}\K(x , x_0)\dx\,.
\end{equation}
Therefore, taking the supremum over $x_0\in\Omega$ of both sides, recalling that $\n$ is increasing,  so that $\sup\limits_{x_0\in\Omega}\n(u(t_0,x_0))=\n\left(\|u(t_0)\|_{\LL^\infty(\Omega)}\right) $, we obtain:
\begin{equation}\label{Upper.PME.Step.2.1}
\n\left(\|u(t_0)\|_{\LL^\infty(\Omega)}\right)
\le \frac{2^{\frac{1}{\mu_0}}}{t_0}\|u(t_0)\|_{\LL^\infty(\Omega)}\sup_{x_0\in\Omega}\int_{\Omega}\K(x , x_0)\dx
\le c_{2,\Omega}\frac{2^{\frac{1}{\mu_0}}}{t_0}\|u(t_0)\|_{\LL^\infty(\Omega)}
\end{equation}
where we have also used the bound \eqref{Lem.Green.est.Upper.I} of Lemma \ref{Lem.Green}, namely  $\sup\limits_{x_0\in\Omega}\int_{\Omega}\K(x , x_0)\dx \le c_{2,\Omega}$\,.\\
Next we recall the Young inequality \eqref{Young.F}\,, namely that for all $a,b\ge 0$ and all $\varepsilon>0$ we have $a\,b\, \le \varepsilon \n_\varepsilon(a)+\varepsilon \nl\left(\frac{b}{\varepsilon}\right)$\,. Applying such inequality with $\varepsilon=1/2$ to \eqref{Upper.PME.Step.2.1} gives:
\[
\n\left(\|u(t_0)\|_{\LL^\infty(\Omega)}\right)
\le c_{2,\Omega}\frac{2^{\frac{1}{\mu_0}}}{t_0}\|u(t_0)\|_{\LL^\infty(\Omega)}\le \frac{1}{2}\n\left(\|u(t_0)\|_{\LL^\infty(\Omega)}\right)
+\frac{1}{2}\nl\left(c_{2,\Omega}\frac{2^{\frac{1}{\mu_0}+1}}{t_0}\right)
\]
Therefore,   we have proved the more general form \eqref{thm.Upper.PME.Absolute.F} of the smoothing effects, namely
\begin{equation}\label{Upper.PME.Step.2.1b}
\n\left(\|u(t_0)\|_{\LL^\infty(\Omega)}\right)
\le \nl\left(\frac{K_1}{t_0}\right)\,,\qquad\mbox{for all $t_0>0$\,, with $K_1=c_{2,\Omega}2^{\frac{1}{\mu_0}+1}\,.$}
\end{equation}


\noindent$\bullet~$\textsc{Step 3. }\textit{Proof of inequality $\eqref{thm.Upper.PME.Absolute}$. }
Now we recall that since we are assuming the hypothesis (N1) on the nonlinearity $\n$, then inequality \eqref{Lem.N1.1} of Lemma \ref{Lem.N1.F} reads
\begin{equation}\label{Upper.PME.Step.3.0}
\n\left(\|u(t_0)\|_{\LL^\infty(\Omega)}\right)\ge (\kb\wedge 1)\n(1)\,\|u(t_0)\|_{\LL^\infty(\Omega)}^{m_j}
\end{equation}
with $j=0$ if $\|u(t_0)\|_{\LL^\infty(\Omega)}\ge 1$ and $j=1$ if $\|u(t_0)\|_{\LL^\infty(\Omega)}\le 1$\,, where $m_i=\frac{1}{1-\mu_i}>1$  (here we really need $\mu_1<1$). On the other hand, inequality \eqref{Lem.N1.1} of Lemma \ref{Lem.N1.F} applied to  $\nl$, reads
\begin{equation}\label{Upper.PME.Step.3.1}
\nl\left(\frac{K_1}{t_0}\right)\le (\ka\vee 1)\,\frac{\nl(r_0)}{r_0^{\frac{m_i}{m_i-1}}}\left(\frac{K_1}{t_0}\right)^{\frac{m_i}{m_i-1}}
\end{equation}
with $i=0$ if $t_0\le K_1/r_0$ and $i=1$ if $t_0\ge K_1/r_0$ \,, where $m_i=\frac{1}{1-\mu_i}>1$.

\noindent Applying the inequalities \eqref{Upper.PME.Step.3.0} and \eqref{Upper.PME.Step.3.1} to \eqref{Upper.PME.Step.2.1b} gives for all $t_0>0$
\begin{equation}\label{Upper.PME.Step.3.2}
(\kb\wedge 1)\n(1)\,\|u(t_0)\|_{\LL^\infty(\Omega)}^{m_j}
\le\n\left(\|u(t_0)\|_{\LL^\infty(\Omega)}\right)
\le \nl\left(\frac{K_1}{t_0}\right)
\le (\ka\vee 1)\,\frac{\nl(r_0)}{r_0^{\frac{m_i}{m_i-1}}}\left(\frac{K_1}{t_0}\right)^{\frac{m_i}{m_i-1}}
\end{equation}
that is
\begin{equation}\label{Upper.PME.Step.3.3}
\|u(t_0)\|_{\LL^\infty(\Omega)}^{m_j}
\le \frac{(\ka\vee 1)}{(\kb\wedge 1)\n(1)\,}\,\frac{\nl(r_0)}{r_0^{\frac{m_i}{m_i-1}}}\left(\frac{K_1}{t_0}\right)^{\frac{m_i}{m_i-1}}
:=\frac{K(r_0, m_i)}{t_0^{\frac{m_i}{m_i-1}}}
\end{equation}
with $i=1$ if $K_1/t_0\le r_0$ and $i=0$ if $K_1/t_0\ge r_0$ and with $j=0$ if $\|u(t_0)\|_{\LL^\infty(\Omega)}\ge 1$ and $j=1$ if $\|u(t_0)\|_{\LL^\infty(\Omega)}\le 1$\,, where $m_i=\frac{1}{1-\mu_i}>1$\,.

The first consequence of \eqref{Upper.PME.Step.3.3} can be obtained by choosing $r_0=K_1$ so that
\begin{equation}\label{Upper.PME.Step.3.4}
\|u(t_0)\|_{\LL^\infty(\Omega)}^{m_j}
\le \frac{(\ka\vee 1)}{(\kb\wedge 1)\n(1)\,}\,\frac{\nl(K_1)}{K_1^{\frac{m_i}{m_i-1}}}\left(\frac{K_1}{t_0}\right)^{\frac{m_i}{m_i-1}}
=\frac{(\ka\vee 1)}{(\kb\wedge 1)\n(1)\,}\,\frac{\nl(K_1)}{t_0^{\frac{m_i}{m_i-1}}}:=\frac{K'_0}{t_0^{\frac{m_i}{m_i-1}}}
\end{equation}
with $i=0$ if $t_0\le 1$ and $i=1$ if $t_0\ge 1$\,. Therefore, there exists a time $\tau_1$ such that $\|u(t_0)\|_{\LL^\infty(\Omega)}\le 1$ for all $t_0\ge \tau_1$\,. On one hand, $\tau_1$ can be zero (for example when $\|u_0\|_{\LL^\infty(\Omega)}\le 1$)\,, therefore it is reasonable to assume that in general $\tau_1$ depends on $u_0$. On the other hand, inequality \eqref{Upper.PME.Step.3.4} shows that for any
\[
t_0\ge \max_{i=0,1} {K'_0}^{\frac{m_i-1}{m_i}}={K'_0}^{\frac{m_1-1}{m_i}}\qquad\mbox{we have\qquad $\|u(t_0)\|_{\LL^\infty(\Omega)}\le 1$}
\]
since it is not restrictive to assume $K_0\ge 1$\,. Therefore we can bound $\tau_1$ from above with a quantity that does not depend on $u_0$, namely
\[
0\le \tau_1(u_0)\le {K'_0}^{(m_1-1)/m_1}:=K_0
\]
Now we fix $r_0=K_1/K_0$\,, so that for $t_0\ge K_0$ we can let $i=j=1$ in inequality \eqref{Upper.PME.Step.3.3}\,, because in this case we have $\|u(t_0)\|_{\LL^\infty(\Omega)}\le 1$ and $K_1/t_0\le r_0$\,, therefore
\begin{equation}\label{Upper.PME.Step.3.3b}
\|u(t_0)\|_{\LL^\infty(\Omega)}^{m_1}
\le \frac{(\ka\vee 1)}{(\kb\wedge 1)\n(1)\,}\,\frac{\nl(K_1/K_0)}{(K_1/K_0)^{\frac{m_1}{m_1-1}}}\left(\frac{K_1}{t_0}\right)^{\frac{m_1}{m_1-1}}
=\frac{{K'_2}^{m_1}}{t^{\frac{m_1}{m_1-1}}}
\end{equation}
Therefore, we have proven inequality \eqref{thm.Upper.PME.Absolute} for $t_0\ge K_0$. Now it remains to prove it for $t_0\le K_0$: first we observe that \eqref{Upper.PME.Step.3.3b} at $t_0=K_0$ gives:
\begin{equation}\label{Upper.PME.Step.3.3c}
\|u(K_0)\|_{\LL^\infty(\Omega)}
\le \left(\frac{(\ka\vee 1)\nl(K_1/K_0)}{(\kb\wedge 1)\n(1)\,}\,\right)^{\frac{1}{m_1}}=K''_2\,.
\end{equation}
then we use the monotonicity inequality \eqref{CP.Monotonicity.3}, namely the function $t\mapsto t^{\frac{1}{m_0-1}}\,u(t,x)$ is nondecreasing in $t>0$ for a.e. $x\in \Omega$\,, from which it follows that for all $0\le t_0\le K_0$ we have
\begin{equation}\label{Upper.PME.Step.2.4}
t_0^{\frac{1}{m_0-1}}\,\|u(t_0)\|_{\LL^\infty(\Omega)}\le K_0^{\frac{1}{m_0-1}}\,\|u(K_0)\|_{\LL^\infty(\Omega)}\le K_0^{\frac{1}{m_0-1}}\,K''_2
\end{equation}
Therefore, we have proven inequality \eqref{thm.Upper.PME.Absolute} for $t_0\le K_0$.

\noindent$\bullet~$\textsc{Step 4. }\textit{Values of the constants. }
\begin{equation}\label{Upper.PME.Step.4}
K_0=\left(\frac{(\ka\vee 1)\nl(K_1)}{(\kb\wedge 1)\n(1)}\right)^{(m_1-1)/m_1}\,,\qquad K_1=c_{2,\Omega}2^{\frac{1}{\mu_0}+1}
\end{equation}
with $c_{2,\Omega}$ is the constant of the bound \eqref{Lem.Green.est.Upper.I} of Lemma \ref{Lem.Green}.
\[\begin{split}
K_2 &= K'_2\,\vee\,\left(K_0^{\frac{1}{m_0-1}}\,K''_2\right)
=\left(\frac{(\ka\vee 1)\nl(K_1/K_0)}{(\kb\wedge 1)\n(1)\,}\,\right)^{\frac{1}{m_1}}K_0^{\frac{m_1-1}{m_1}}\mbox{\,.\qed}
\end{split}\]

\noindent\textbf{Remark. }The precise boundary behaviour will be studied in the forthcoming paper \cite{BV-ppr2-2}.


\section{Smoothing effects}\label{sect.smooth}

In this Section we prove three formulas for the bounds called ``smoothing effects'',  where the $L^\infty$ norm of the solution at time $t>0$ is estimated in terms of possibly weighted integral norms of the initial data, or even of the solution at the same time or at a previous times (instantaneous or backward smoothing effect). The instantaneous and backward smoothing effects bounds are new also in the well studied case $s=1$. We recall the exponents for $\gamma\in [0,1]$:
\[
\vartheta_{i,\gamma}=\frac{1}{2s+(N+\gamma)(m_i-1)}\qquad\mbox{with }\qquad m_i=\frac{1}{1-\mu_i}>1
\]
In the next results, $\gamma$ is the exponent appearing in assumption $(K2)$.
\begin{thm}[Weighted $\LL^1-\LL^\infty$ smoothing effect]\label{thm.Upper.1.PME} As a consequence of $(K2)$ hypothesis, there exists a constant $K_6>0$ such that the following estimates hold true.
\begin{equation}\label{thm.Upper.PME.Smoothing.F}
\n\big(\|u(t)\|_{\LL^\infty(\Omega)}\big)
\le K_6\, \frac{\| u(t_0)\|_{\LL^1_{\p}(\Omega)}^{2s m_i\vartheta_{i,\gamma}}}{t^{m_i(N+\gamma)\vartheta_{i,\gamma}}}\,,
\qquad\mbox{for all }0\le t_0\le t\,,
\end{equation}
with $i=1$ if $t\ge \| u(t_0)\|_{\LL^1_{\p}(\Omega)}^{\frac{2s}{N+\gamma}}$ and $i=0$
if $t\le \| u(t_0)\|_{\LL^1_{\p}(\Omega)}^{\frac{2s}{N+\gamma}}$\,.
\end{thm}
The estimate is a consequence of the $(K2)$ bounds for the Green function. The constant $K_6 >0$ depends only on $d, m_i, s,\gamma, F$ and $\Omega$\,, and have an explicit form given in the proof. $K_1>0$ is the universal constant of Theorem $\ref{thm.Upper.PME}$\,.
\begin{cor}\label{cor.smoothing.0}Under the weaker assumption $(K1)$ instead of $(K2)$, the results of Theorem $\ref{thm.Upper.1.PME}$ hold true with $\gamma=0$ and replacing $\|\cdot\|_{\LL^1_{\p}(\Omega)}$ with $\|\cdot\|_{\LL^1(\Omega)}$\,.
\end{cor}
\noindent {\bf Proof of Theorem \ref{thm.Upper.1.PME}.} Let $x_0\in \Omega$ and consider $B_r(x_0)$ with $r>0$ to be fixed later. Define the set $\Omega_r=\Omega\setminus\left(B_r(x_0)\cap \Omega\right)$ so that $\Omega\subseteq B_r(x_0)\cup \Omega_r$. Notice that the ball $B_r(x_0)$ need not to be included in $\Omega$\,. Then it is clear that $\forall x\in\Omega_r$ we have $|x-x_0|\ge r$\,. The starting point of the argument is the fundamental upper estimate \eqref{Upper.PME.Step.1.2} of Step 1 of the proof of Theorem \ref{thm.Upper.PME}
combined with assumption (K2). To be precise, we fix $t_0>0$, and use these results to get
\begin{equation*}\begin{split}
\n\left(u(t_0,x_0)\right) &\le \frac{2^{\frac{1}{\mu_0}}}{t_0}\int_{\Omega}u(t_0,x)\K(x , x_0)\dx\\
&\le \frac{2^{\frac{1}{\mu_0}}}{t_0}\left[\int_{B_r(x_0)}u(t_0,x)\K(x , x_0)\dx+\int_{\Omega_r}u(t_0,x)\K(x , x_0)\dx\right]\\
~_{(a)}&\le c_{1,\Omega}\frac{2^{\frac{1}{\mu_0}}}{t_0}
    \left[\|u(t_0)\|_{\LL^\infty(\Omega)}\int_{B_r(x_0)}\frac{1}{|x-x_0|^{N-2s}}\dx
        +\int_{\Omega_r}\frac{u(t_0,x)\p(x)}{|x-x_0|^{N-2s+\gamma}}\dx\right]\\
~_{(b)}&\le \frac{1}{2}\n\big(\|u(t_0)\|_{\LL^\infty(\Omega)}\big)
        + \frac{1}{2}\nl\left( \frac{\omega_d\,c_{1,\Omega}2^{\frac{1}{\mu_0}}}{s} \frac{r^{2s}}{t_0}\right)
        + \frac{c_{1,\Omega}2^{\frac{1}{\mu_0}}}{t_0\,r^{N-2s+\gamma}}\int_{\Omega}u(t_0,x)\p(x)\dx \\
\end{split}
\end{equation*}
where in $(a)$ we have used  the Green function estimates (K1)\,, namely $\K(x,x_0)\le c_{1,\Omega}|x-x_0|^{-(N-2s)}$ on the ball $B_r(x_0)$\,, while on $\Omega_r$ we have used the Green function estimates (K2)\,,  namely $$\K(x,x_0)\le c_{1,\Omega}\p(x)|x-x_0|^{-(N-2s+\gamma)}\,.
$$
In $(b)$  we have used the Young inequality \eqref{Young.F}, namely
\[
ab\le \frac{1}{2}\n(a)+\frac{1}{2}\nl(2b)\qquad\mbox{for all }a,b\ge 0\,.
\]
Taking the supremum over $x_0\in \Omega$ in the above inequality gives for all $r>0$:
\begin{equation}\label{WSME.Step2.1}
\n\big(\|u(t_0)\|_{\LL^\infty(\Omega)}\big)\le \nl\left( \frac{\omega_d\,c_{1,\Omega}2^{\frac{1}{\mu_0}}}{s} \frac{r^{2s}}{t_0}\right)
        + \frac{c_{1,\Omega}2^{\frac{1}{\mu_0}+1}}{t_0\,r^{N-2s+\gamma}}\|u(t_0)\|_{\LL^1_{\p}(\Omega)}
\end{equation}
since we recall that $\sup\limits_{x_0\in \Omega}\n(u(t,x_0))=\n\big(\|u(t_0)\|_{\LL^\infty(\Omega)}\big)$, because $F$ is nondecreasing.

\medskip

\noindent (ii) In order to get a simpler formula we optimize in free parameter $r$. We want to eliminate the first term of the right-hand side. For that, and in view of the behaviour of $\nl$,  we choose a $0\le \tau_0\le t_0$ and let  $i=1$ if $t_0\ge \| u(\tau_0)\|_{\LL^1_{\p}(\Omega)}^{\frac{2s}{N+\gamma}}$ and $i=0$ if $t_0\le \| u(\tau_0)\|_{\LL^1_{\p}(\Omega)}^{\frac{2s}{N+\gamma}}$\,. Put then
\[
r=\left(\frac{s}{\omega_d\,c_{1,\Omega}2^{\frac{1}{\mu_0}}}\right)^{\frac{1}{2s}}
    \left(\tau_0\,\|u(\tau_0)\|_{\LL^1_{\p}(\Omega)}^{m_i-1}\right)^{\vartheta_{i,\gamma}}
\qquad\mbox{with}\qquad \vartheta_{i,\gamma}=\frac{1}{2s-(N+\gamma)(m_i-1)}\,.
\]
The above choice of $r$ implies the following estimates for the two terms in the right-hand side of inequality \eqref{WSME.Step2.1}. As the first terms is concerned, we have:
\[
\frac{\omega_d\,c_{1,\Omega}2^{\frac{1}{\mu_0}}}{s}\frac{r^{2s}}{t_0}=\frac{\tau_0^{2s\vartheta_{i,\gamma}}}{t_0}
    \,\|u(\tau_0)\|_{\LL^1_{\p}(\Omega)}^{2s(m_i-1)\vartheta_{i,\gamma}}
\le \left[\frac{\| u(\tau_0)\|_{\LL^1_{\p}(\Omega)}^{2s\vartheta_{i,\gamma}}}{t_0^{(N+\gamma)\vartheta_{i,\gamma}}}\right]^{m_i-1}
=: H_i^{m_i-1}
\]
where in the last step we have used the absolute upper bound \eqref{thm.Upper.PME.Absolute}.
Indeed, this implies that the right-hand side of \eqref{WSME.Step2.1} becomes
\[
\nl\left( \frac{\omega_d\,c_{1,\Omega}2^{\frac{1}{\mu_0}}}{s} \frac{r^{2s}}{t_0}\right)
\le \nl(H_i^{m_i-1})\le (\overline{k}\vee 1)\, \nl(1)\, H_i^{m_i}
\]
where in the second step we have used Lemma \ref{Lem.N1.F}, which implies that $\nl(X)\le \overline{k}\vee 1\, \nl(1) X^{m_i/(m_i-1)}$ with $i=1$ if $X\in [0,1]$ and $i=0$ if $X\ge 1$\,. This explains our choice of $i$. As a consequence we obtain
\[
\nl\left( \frac{\omega_d\,c_{1,\Omega}2^{\frac{1}{\mu_0}}}{s} \frac{r^{2s}}{t_0}\right)
\le K'_6\left[\frac{\|u(\tau_0)\|_{\LL^1_{\p}(\Omega)}^{2s\vartheta_{i,\gamma}}}
    {t_0^{(N+\gamma)\vartheta_{i,\gamma}}}\right]^{m_i}
\]
with  $i=1$ if $t_0\ge \| u(\tau_0)\|_{\LL^1_{\p}(\Omega)}^{\frac{2s}{N+\gamma}}$ and $i=0$ if $t_0\le \| u(\tau_0)\|_{\LL^1_{\p}(\Omega)}^{\frac{2s}{N+\gamma}}$\,.

\noindent We now estimate the second term in the right-hand side of \eqref{WSME.Step2.1} as follows
\[\begin{split}
\frac{c_{1,\Omega}2^{\frac{1}{\mu_0}+1}}{t_0\,r^{N-2s+\gamma}}\|u(t_0)\|_{\LL^1_{\p}(\Omega)}
&\le \frac{c_{1,\Omega}2^{\frac{1}{\mu_0}+1}}{t_0\, \tau_0^{(N-2s+\gamma)\vartheta_{i,\gamma}}}
\left(\frac{\omega_d\,c_{1,\Omega}2^{\frac{1}{\mu_0}}}{s}\right)^{\frac{N-2s+\gamma}{2s}}
\frac{\|u(t_0)\|_{\LL^1_{\p}(\Omega)}}{\|u(\tau_0)\|_{\LL^1_{\p}(\Omega)}^{(N-2s+\gamma)(m_i-1)\vartheta_{i,\gamma}}}
\end{split}\]
Now we use the ``quasi'' monotonicity of the $\LL^1_{\p}(\Omega)$-norm, namely that $\|u(t_0)\|_{\LL^1_{\p}(\Omega)}\le C_{\Omega,\gamma}\|u(\tau_0)\|_{\LL^1_{\p}(\Omega)}$ for $0\le \tau_0\le t_0$, see formula \eqref{monotonicity.dist} of Proposition \ref{thm.L1weight.contr.psi}. We obtain\normalcolor
\[\begin{split}
\frac{c_{1,\Omega}2^{\frac{1}{\mu_0}+1}}{t_0\,r^{N-2s+\gamma}}\|u(t_0)\|_{\LL^1_{\p}(\Omega)}
&\le  K_6 \frac{\|u(\tau_0)\|_{\LL^1_{\p}(\Omega)}^{1-(N-2s+\gamma)(m_i-1)\vartheta_{i,\gamma}}}{t_0^{1+(N-2s+\gamma)\vartheta_{i,\gamma}}}
=K_6 \left[\frac{\|u(\tau_0)\|_{\LL^1_{\p}(\Omega)}^{2s\vartheta_{i,\gamma}}}{t_0^{(N+\gamma)\vartheta_{i,\gamma}}}\right]^{m_i}
=K_6 H_i^{m_i}
\end{split}
\]
Finally, summing up the two estimates we have obtained that inequality \eqref{WSME.Step2.1} becomes
\begin{equation}\label{WSME.Step2.2}
\n\big(\|u(t_0)\|_{\LL^\infty(\Omega)}\big)\le K'_6 H_i^{m_i}
        + K_6\, H_i^{m_i}\,,
\quad\mbox{with}\quad
H_i= \frac{\| u(\tau_0)\|_{\LL^1_{\p}(\Omega)}^{2s\vartheta_{i,\gamma}}}{t_0^{(N+\gamma)\vartheta_{i,\gamma}}}
\end{equation}
with $i=1$ if $t_0\ge \| u(\tau_0)\|_{\LL^1_{\p}(\Omega)}^{\frac{2s}{N+\gamma}}$ and $i=0$ if $t_0\le \| u(\tau_0)\|_{\LL^1_{\p}(\Omega)}^{\frac{2s}{N+\gamma}}$\,. The proof of the smoothing effect \eqref{thm.Upper.PME.Smoothing.F} is concluded after changing the letters: $t_0$ into $t$ and $\tau_0$ into $t_0$.\qed

\noindent {\bf Proof of Corollary \ref{cor.smoothing.0}. }The proof is a consequence of hypothesis (K1), and it has been given in \cite{BV-PPR1} for the case $\n(u)=u^m$ and $\gamma=1$. We also remark that the proof is the same as above, just letting $\gamma=0$ and using hypothesis (K1) instead of (K2)\,.\qed

\noindent\textbf{Remarks. }(i) The power in the smoothing effects changes for small and large times, more precisely it depends on  $\|u(t_0)\|_{\LL^\infty(\Omega)}$. We may choose $t_0=0$ to have an explicit dependence on the data. An important point is that the power of time is always less than one, more specifically, $(d+\gamma)(m_i-1)\vartheta_{i,\gamma}<1$ for $i=0,1$\,; this fact has important consequences, for instance is crucial in the proof of $\LL^1_{\p}$ weighted estimates of Proposition \ref{thm.L1weight.contr}, which are essential to prove existence of weak dual solutions, Theorem \ref{thm.L1weight.exist}; they are crucial also in the proof of the lower bounds of the forthcoming paper \cite{BV-ppr2-2}.\\
(ii) The weighted smoothing effects are new to our knowledge also when $s=1$. Moreover, they apply to a class of nonnegative initial data $\LL^1_{\p}(\Omega) $ which is strictly larger than $\LL^1(\Omega)$.\\
(iii) Another novelty is represented by the fact that the smoothing effect occurs at the same time; this is new also when $s=1$\,.\\
(iv) In the case when $\n(u)=u^m$\,, we recover the sharp results of our previous paper \cite{BV-PPR1}\,.\\
(v) The smoothing effect of Theorem has an intrinsic form, that can be made more explicit by separating the result for small and large times as in the following corollary.\normalcolor
\begin{cor}\label{cor.smoothing.1}There exists a constant $K_7>0$ such that the following estimates hold true.\\
\noindent \textsc{Weighted $\LL^1-\LL^\infty$ smoothing effect for small times: }
\begin{equation}\label{thm.Upper.PME.Smoothing.small}
\|u(t)\|_{\LL^\infty(\Omega)}
\le K_7\, \frac{\| u(t_0)\|_{\LL^1_{\p}(\Omega)}^{2s \vartheta_{0,\gamma}}}{t^{(N+\gamma)\vartheta_{0,\gamma}}}\,,
\qquad\mbox{for all $0\le t_0\le t \le \| u(t_0)\|_{\LL^1_{\p}(\Omega)}^{\frac{2s}{N+\gamma}}\,.$}
\end{equation}
\noindent \textsc{Weighted $\LL^1-\LL^\infty$ smoothing effect for large times: }
\begin{equation}\label{thm.Upper.PME.Smoothing.large}
\|u(t)\|_{\LL^\infty(\Omega)}
\le K_7\, \frac{\| u(t_0)\|_{\LL^1_{\p}(\Omega)}^{2s \vartheta_{1,\gamma}}}{t^{(d+\gamma)\vartheta_{1,\gamma}}}\,,
\qquad\mbox{for all $t\ge  \| u(t_0)\|_{\LL^1_{\p}(\Omega)}^{\frac{2s}{d+\gamma}}\,.$}
\end{equation}
Moreover, the condition  $t\ge  \| u(t_0)\|_{\LL^1_{\p}(\Omega)}^{\frac{2s}{d+\gamma}}$\,, is implied by $t\ge \left(K_1\,\|\p\|_{\LL^1(\Omega)} \right)^{\vartheta_{1,\gamma}(m_1-1)}$\,.
\end{cor}
The constant $K_7 >0$ depends only on $d, m_i, s,\gamma, F$ and $\Omega$\,, and have an explicit form given in the proof. $K_1>0$ is the universal constant of Theorem $\ref{thm.Upper.PME}$\,.

 As it happened for Corollary \ref{cor.smoothing.0}, if we only assume (K1), then the result of Corollary \ref{cor.smoothing.1} hold true also for $\gamma=0$, just by replacing $\|\cdot\|_{\LL^1_{\p}(\Omega)}$ with $\|\cdot\|_{\LL^1(\Omega)}$\,.

\noindent {\sl Proof.~}We split the proof in two steps.

\noindent$\bullet~$\textit{Weighted smoothing effect for large times. }We first prove the instantaneous smoothing effects, namely $t_0=\tau_0$. When $t_0\ge\|u(t_0)\|_{\LL^1_{\p}(\Omega)}^{\frac{2s}{N+\gamma}}$\,, then $i=1$, and we know that
\[
\|u(t_0)\|_{\LL^\infty(\Omega)}\ge \frac{\|u(t_0)\|_{\LL^1_{\p}(\Omega)}}{\|\p\|_{\LL^1(\Omega)}}:=r_0
\ge \frac{t_0^{\frac{N+\gamma}{2s}}}{\|\p\|_{\LL^1(\Omega)}}>0
\]
Therefore, using Lemma \ref{Lem.N1.F}\,, we get that $F(X)\ge \kb \n(r_0)(X/r_0)^{m_1}\ge \kb F(1) X^{m_1}$, since
$\lim\limits_{r\to 0}\frac{\n(r)}{r^{m_1}}\ge \n(1)>0$. As a consequence, we get:
\[
\kb\n(1)\|u(t_0)\|_{\LL^\infty(\Omega)}^{m_1}\le \n\big(\|u(t_0)\|_{\LL^\infty(\Omega)}\big)
\le(1+K_6)\,\frac{\| u(t_0)\|_{\LL^1_{\p}(\Omega)}^{2s m_1\vartheta_{1,\gamma}}}{t_0^{m_1(N+\gamma)\vartheta_{1,\gamma}}} \qquad\mbox{for all }t_0\ge\|u(t_0)\|_{\LL^1_{\p}(\Omega)}^{\frac{2s}{N+\gamma}}\,.
\]
We have proven \eqref{thm.Upper.PME.Smoothing.large} with $K_7=(1+K_6)/[\kb\,\n(1)]$ when $\tau_0=t_0$.

When $0\le \tau_0< t_0$ we obtain \eqref{thm.Upper.PME.Smoothing.large} using the ``quasi'' monotonicity of the $\LL^1_{\p}(\Omega)$-norm, namely that $\|u(t_0)\|_{\LL^1_{\p}(\Omega)}\le C_{\Omega,\gamma}\|u(\tau_0)\|_{\LL^1_{\p}(\Omega)}$ for $0\le \tau_0\le t_0$, see formula \eqref{monotonicity.dist} of Proposition \ref{thm.L1weight.contr.psi}.

Finally, we can provide a more explicit condition. Indeed, a close inspection of the proof reveals that the true condition for bounds \eqref{thm.Upper.PME.Smoothing.large} to hold is $t\ge  \| u(t)\|_{\LL^1_{\p}(\Omega)}^{\frac{2s}{N+\gamma}}\,.$ Thanks to the absolute bounds of Theorem \ref{thm.Upper.PME}, we then get that
\[
\| u(t)\|_{\LL^1_{\p}(\Omega)}\le \frac{K_1\,\|\p\|_{\LL^1(\Omega)} }{t^{\frac{1}{m_1-1}}}\le t^{\frac{N+\gamma}{2s}}
\qquad\mbox{is implied by}\qquad
t^{\frac{1}{\vartheta_{1,\gamma}(m_1-1)}}=t^{\frac{N+\gamma}{2s}+\frac{1}{m_1-1}}\ge \left(K_1\,\|\p\|_{\LL^1(\Omega)} \right),
\]
therefore, the bounds \eqref{thm.Upper.PME.Smoothing.large.backw} hold for any $t\ge \left(K_1\,\|\p\|_{\LL^1(\Omega)} \right)^{\vartheta_{1,\gamma}(m_1-1)}$, where $K_1>0$ is the universal constant given in Theorem \ref{thm.Upper.PME}\,.

\noindent$\bullet~$\textit{Weighted smoothing effect for small times. }When $t_0\le \| u(\tau_0)\|_{\LL^1_{\p}(\Omega)}^{\frac{2s}{N+\gamma}}$ then $i=0$\,, and $U_\gamma(\tau_0, t_0)\ge 1$. We split two cases, namely  $\|u(t_0)\|_{\LL^\infty(\Omega)}\le U_\gamma(\tau_0, t_0)$ and $\|u(t_0)\|_{\LL^\infty(\Omega)}\ge U_\gamma(\tau_0, t_0)$. In the first case, there is nothing to prove, indeed
\[
\|u(t_0)\|_{\LL^\infty(\Omega)}\le U_\gamma(\tau_0, t_0)\qquad\mbox{means}\qquad
\|u(t_0)\|_{\LL^\infty(\Omega)}\le \frac{\| u(\tau_0)\|_{\LL^1_{\p}(\Omega)}^{2s\vartheta_{0,\gamma}}}{t_0^{(N+\gamma)\vartheta_{0,\gamma}}}\,,
\]
which is \eqref{thm.Upper.PME.Smoothing.small} with $K_7=1$\,. In the second case, we take $\|u(t_0)\|_{\LL^\infty(\Omega)}\ge U_\gamma(\tau_0, t_0)\ge 1$. By Lemma \ref{Lem.N1.F}, $\n(X)\ge \n(1)X^{m_0}$ when $X\ge 1$\,, so that
\[
\n\big(\|u(t_0)\|_{\LL^\infty(\Omega)}\big)
\ge F(1)\|u(t_0)\|_{\LL^\infty(\Omega)}^{m_0}\,.
\]
We have proven \eqref{thm.Upper.PME.Smoothing.small} with $K_7=(1+K_6)/F(1)$.
Finally we can put $K_7:= 1\vee (1+K_6)/(\kb F(1)\wedge 1)$.\qed

\noindent As a further corollary of Theorem \ref{thm.Upper.1.PME}, we get the following reverse in time smoothing effects.
\begin{cor}[Backward Smoothing effects]\label{thm.Upper.Backward.PME} As a consequence of $(K2)$ hypothesis, there exists a constant $K_7>0$ such that the following estimates hold true for all $t,h>0$:

\noindent\textsc{Backward $\LL^1-\LL^\infty$-weighted smoothing effect for small times: }
\begin{equation}\label{thm.Upper.PME.Smoothing.small.backw}
\|u(t)\|_{\LL^\infty(\Omega)}
\le 2 K_7\, \left(1\vee\frac{h}{t}\right)^{\frac{2s \vartheta_{0,\gamma}}{m_0-1}} \frac{\| u(t+h)\|_{\LL^1_{\p}(\Omega)}^{2s \vartheta_{0,\gamma}}}{t^{(N+\gamma)\vartheta_{0,\gamma}}}\,,
\qquad\mbox{for all $0\le t \le \| u(t)\|_{\LL^1_{\p}(\Omega)}^{\frac{2s}{N+\gamma}}\,.$}
\end{equation}

\noindent\textsc{Backward $\LL^1-\LL^\infty$-weighted smoothing effect for large times: }
\begin{equation}\label{thm.Upper.PME.Smoothing.large.backw}
\|u(t)\|_{\LL^\infty(\Omega)}
\le 2 K_7\,\left(1\vee\frac{h}{t}\right)^{\frac{2s \vartheta_{1,\gamma}}{m_1-1}}
 \frac{\| u(t+h)\|_{\LL^1_{\p}(\Omega)}^{2s \vartheta_{1,\gamma}}}{t^{(N+\gamma)\vartheta_{1,\gamma}}}\,,
\qquad\mbox{for all $t \ge \| u(t)\|_{\LL^1_{\p}(\Omega)}^{\frac{2s}{N+\gamma}}\,.$}
\end{equation}
Moreover, the condition  $t\ge  \| u(t)\|_{\LL^1_{\p}(\Omega)}^{\frac{2s}{N+\gamma}}$\,, is implied by $t\ge \left(K_1\,\|\p\|_{\LL^1(\Omega)} \right)^{\vartheta_{1,\gamma}(m_1-1)}$\,.\end{cor}

\noindent\textbf{Remarks. } The constant $K_7>0$ do not depend on $u$ and is given in Theorem $\ref{thm.Upper.1.PME}$. $K_1>0$ is the universal constant of Theorem $\ref{thm.Upper.PME}$\,. Moreover, for $\gamma=0$ we obtain the smoothing effects for $u_0\in \LL^1(\Omega)$, as consequence only of $(K1)$ bounds, and the above formulas hold replacing $\|\cdot\|_{\LL^1_{\p}(\Omega)}$ with $\|\cdot\|_{\LL^1(\Omega)}$\,.

\noindent {\sl Proof of  Corollary \ref{thm.Upper.PME.Smoothing.small.backw}.~}We use the monotonicity estimates given by Theorem \ref{thm.CP-JFA}, namely \eqref{CP.Monotonicity.3}, which imply that the function $t\mapsto t^{\frac{1-\mu_0}{\mu_0}}\,u(t,x)=t^{\frac{1}{m_0-1}}\,u(t,x)$ is nondecreasing in $t>0$ for a.e. $x\in \Omega$\,. Indeed we have that for all $h\ge 0$ and $t\ge 0$:
\[
\| u(t)\|_{\LL^1_{\p}(\Omega)}^{2s \vartheta_{1,\gamma}}
\le \left(\frac{t+h}{t}\right)^{\frac{2s \vartheta_{1,\gamma}}{m_0-1}}\| u(t+h)\|_{\LL^1_{\p}(\Omega)}^{2s \vartheta_{1,\gamma}}
\le 2\left(1\vee\frac{h}{t}\right)^{\frac{2s \vartheta_{1,\gamma}}{m_0-1}}\| u(t+h)\|_{\LL^1_{\p}(\Omega)}^{2s \vartheta_{1,\gamma}}
\]
The proof is concluded once we apply the above inequality to the smoothing effects of Theorem \ref{thm.Upper.1.PME} with $t_0=t$\,.\qed

\medskip

\noindent\textbf{Remark. }As a concluding remark for this section about upper estimates, we shall observe that in order to get estimates for $\n(u) \in \LL^1\left((0,\infty):\LL^1_{\p}(\Omega)\right)$, a space required by our concept of solution, we need both the absolute bound and the smoothing effect of Theorems \eqref{thm.Upper.PME} and \eqref{thm.Upper.1.PME} respectively. Let $t_0=\| u(\tau_0)\|_{\LL^1_{\p}(\Omega)}^{\frac{2s}{N+\gamma}}$
\begin{equation}\label{remark.upper.final}\begin{split}
0\le & \int_0^{\infty}\int_{\Omega} \n(u(t,x))\p(x)\dx\dt
\le \int_0^{K_0}\int_{\Omega} \n(u(t,x))\p(x)\dx\dt +\int_{K_0}^{\infty}\int_{\Omega} \n(u(t,x))\p(x)\dx\dt \\
&\le \int_0^{t_0}\frac{K_6^{m_0-1}}{t^{(N+\gamma)(m_0-1)\vartheta_{1,\gamma}}}\|u_0\|_{\LL^1_{\p}(\Omega)}^{2s(m-1)\vartheta_{1,\gamma}+1}\dt +\int_{t_0}^{\infty}\frac{K_1^{m_1}\|\p\|_{\LL^1(\Omega)}}{t^{\frac{m_i}{m_i-1}}}\dt <+\infty\\
\end{split}\end{equation}
the first integral is finite since $(N+\gamma)(m_0-1)\vartheta_{1,\gamma}<1$ and the second since $m_i/(m_i-1)>1$\,. We have also used the ``quasi'' monotonicity of the $\LL^1_{\p}(\Omega)$-norm, see formula \eqref{monotonicity.dist} of Proposition \ref{thm.L1weight.contr.psi}.

\section{Mild solutions vs weak dual solutions}\label{Sect.Various.Sols}

In this section we are going to prove that the class $\mathcal{S}_p$ of weak dual solutions for which our estimates hold (at the first step of the proofs) contains all mild solutions with initial data $u_0\in \LL^p(\Omega)$. To this end we prove two propositions which have their own interest, and when joined together they prove our claim. The first proposition proves that the $\LL^p$-norm of the mild solution constructed in Theorem \ref{thm.CP-JFA} is non-increasing in time, therefore $u(t)\in\LL^p(\Omega)$ when $u_0$ does. The second proposition shows that mild solutions are indeed weak dual solutions according to Definition \ref{Def.Very.Weak.Sol.Dual} when $u_0\in \LL^p(\Omega)$, with $p\in [1,\infty]$\,.

\begin{prop}[Semigroup solutions belong to $\LL^p(\Omega)$ when $u_0$ does]\label{Prop.Lp.Mild}
Let $u$ be the unique semigroup (mild) solution corresponding to the initial datum $u_0\in \LL^p(\Omega)$ with $p\ge 1$\,, as in Theorem $\ref{thm.CP-JFA}$. Then $u(t)\in \LL^p(\Omega)$ for all $t>0$, more precisely $\|u(t)\|_{\LL^p(\Omega)}\le \|u_0\|_{\LL^p(\Omega)}$\,.
\end{prop}

\noindent {\sl Proof.~}If $p=1$ the result is true since the semigroup is contractive on $\LL^1$. From now on let us fix $1<p<\infty$ and assume $u_0\in \LL^p(\Omega)$. Indeed, it is not restrictive to assume $u_0\in  \LL^\infty(\Omega)\subset \LL^p(\Omega)$\,; consider an approximating sequence $u_{0,j}\in   \LL^\infty(\Omega)$ that converges strongly to $u_0\in \LL^p(\Omega)$ for which $\|u_j(t)\|_{\LL^p(\Omega)}\le \|u_{0,j}\|_{\LL^p(\Omega)}$; letting $j\to\infty$ we obtain $\|u(t)\|_{\LL^p(\Omega)}\le \|u_0\|_{\LL^p(\Omega)}$\,, using the lower semicontinuity of the norm. Once the result is true for all $p>1$ then the case $p=\infty$ follows by taking the limit $p\to\infty$\,.\normalcolor

Assume that $u_0\in \LL^\infty(\Omega)$ and $1<p<\infty$. The semigroup (mild) solutions of Theorem \ref{thm.CP-JFA} are obtained via an implicit time discretization method (or via the Crandall-Liggett theorem), as follows. Consider the following partition of $[0,T]$
\[
t_k= \frac{k}{n}T\,,\qquad\mbox{for any }0\le k\le n\qquad\mbox{(recall that $t_0=0$ and $t_n=T$)}
\]
and let finally $h=t_{k+1}-t_k=T/n$. For any $t\in (0,T)$\,, the (unique) semigroup solution $u(t,\cdot)$ is obtained as the limit in $\LL^1(\Omega)$ of the solutions $u_{k+1}(\cdot)=u(t_{k+1},\cdot)$ which solve the following elliptic equation ($u_k$ is the datum, is given by the previous iterative step)
\begin{equation}\label{1.1.Lp.Stab}
h\A\n(u_{k+1})+u_{k+1}=u_k\qquad\mbox{or equivalently}\qquad\frac{u_{k+1}-u_k}{h}=-\A\n(u_{k+1})\,.
\end{equation}
Moreover, since we are assuming $u_0\in \LL^\infty(\Omega)$, all the approximations $u_k\in \LL^\infty(\Omega)$.   This is a consequence of the rather standard elliptic estimates, see the proofs of \cite{CP-JFA} for more details.\normalcolor

We multiply the above equation by $u_{k+1}^{p-1}=u^{p-1}(t_{k+1},\cdot)$, we integrate and we sum over $k$: we obtain by rearranging the sums
\begin{equation}\label{1.2.Lp.Stab}\begin{split}
- \int_\Omega u_{k+1}^{p-1}\A\n(u_{k+1})\dx
&= \int_\Omega (u_{k+1}-u_k)u_{k+1}^{p-1}\dx\\
\end{split}
\end{equation}
We first observe that as a consequence of Assumption (A2'), we have that
\[
\int_\Omega u_{k+1}^{p-1}\A\n(u_{k+1})\dx = \int_\Omega \n\left(w^{\frac{1}{p-1}}\right)\A w\dx= \int_\Omega \beta(w) \A w\dx\ge 0
\]
notice that we have let $w=u_{k+1}^{p-1}$, and the above inequality is true by (A2') since  the function $w\mapsto \beta(w)=\n\left(w^{\frac{1}{p-1}}\right)$ is  a  monotone and continuous real function with $\beta(0)=0$ for all $p>1$, since both parts are increasing. Since $w$ and $\A w \in \LL^\infty(\Omega)$ we can apply (A2'). Therefore,
\begin{equation}
\int_{\Omega}u_{k+1}^p\dx\le \int_{\Omega}u_k\, u_{k+1}^{p-1}\dx
\end{equation}
so that by H\"older inequality, we get $\|u_{k+1}\|_{\LL^p(\Omega)}\le \|u_k\|_{\LL^p(\Omega)}$ for every $k$. Letting $k\to\infty$ we obtain the monotonicity of the $\LL^p$ norm for mild solutions.\qed

\medskip

\begin{prop}[Semigroup solutions are  weak dual solutions]\label{Prop.Mild-WDS}
Let $u$ be the unique semigroup (mild) solution corresponding to the initial datum $u_0\in \LL^1(\Omega)$\,, as in Theorem $\ref{thm.CP-JFA}$. Then $u$ is a weak dual solution in the sense of Definition $\ref{Def.Very.Weak.Sol.Dual}$\,.
\end{prop}
\noindent {\sl Proof.~}Assume that $u_0\in \LL^\infty(\Omega)$. The semigroup (mild) solutions of Theorem \ref{thm.CP-JFA} are obtained via an implicit time discretization method (or via the Crandall-Liggett theorem), as follows. Consider the following partition of $[0,T]$
\[
t_k= \frac{k}{n}T\,,\qquad\mbox{for any }0\le k\le n\qquad\mbox{(recall that $t_0=0$ and $t_n=T$)}
\]
and let finally $h=t_{k+1}-t_k=T/n$. For any $t\in (0,T)$\,, the (unique) semigroup solution $u(t,\cdot)$ is obtained as the limit in $\LL^1(\Omega)$ of the solutions $u_{k+1}(\cdot)=u(t_{k+1},\cdot)$ which solve the following elliptic equation ($u_k$ is the datum, is given by the previous iterative step)
\begin{equation}\label{1.1.Mild-WDS}
h\A\n(u_{k+1})+u_{k+1}=u_k\qquad\mbox{or equivalently}\qquad\frac{u_{k+1}-u_k}{h}=-\A\n(u_{k+1})\,.
\end{equation}\normalcolor
Moreover, since we are assuming $u_0\in \LL^\infty(\Omega)$, all the approximations $u_k\in \LL^\infty(\Omega)$.  This is a consequence of the rather standard elliptic estimates, see the proofs of \cite{CP-JFA} for more details.

Therefore we can apply $\AI$ to both members, to obtain
\begin{equation}\label{1.2.Mild-WDS}
h\n(u_{k+1})+\AI u_{k+1}=\AI u_k\qquad\mbox{or equivalently}\qquad\frac{\AI u_{k+1}-\AI u_k}{h}=-\n(u_{k+1})\,.
\end{equation}\normalcolor
Multiply both members of equality \eqref{1.2.Mild-WDS} by $\psi_k=\psi(t_k,\cdot)$ where $\psi$ is an admissible test function in the sense of Definition \ref{Def.Very.Weak.Sol.Dual}, namely $\psi/\p\in C^1_c((0,T): \LL^\infty(\Omega))$; an integration over $\Omega$ and summation over $k$ gives
\begin{equation}\label{1.3.Mild-WDS}
\sum_{k=0}^{n-1}\int_{\Omega}(\AI u_{k+1}- \AI u_k)\psi_k\dx
=-\sum_{k=0}^{n-1}h\int_{\Omega}\n(u_{k+1})\psi_k\dx\,.
\end{equation}
Notice that it is not restrictive to assume (at least for $n$ large)
\begin{equation}\label{psi.cpt}
\psi_1=\psi(t_1,\cdot)=\psi(T/n,\cdot)=0\qquad\mbox{and}\qquad \psi(T,\cdot)=0\,.
\end{equation}
since $\psi$ is compactly supported (in time) in $(0,T)$\,. Next we observe that:
\begin{equation}\label{1.4.Mild-WDS}
\sum_{k=0}^{n-1}\int_{\Omega}(\AI u_{k+1}- \AI u_k)\psi_k\dx
= \int_{\Omega}\left(\psi_n \AI u_n -\psi_1\AI u_0\right)\dx -\sum_{k=1}^{n-1}\int_{\Omega}(\psi_{k+1}-\psi_k)\AI u_k\dx\,,
\end{equation}
combining this with \eqref{1.3.Mild-WDS} we get:
\begin{equation}\label{1.5.Mild-WDS}
\int_{\Omega}\left(\psi_n \AI u_n -\psi_1\AI u_0\right)\dx - \sum_{k=1}^{n-1}\int_{\Omega}(\psi_{k+1}-\psi_k)\AI u_k\dx
=-\sum_{k=0}^{n-1}h\int_{\Omega}\n(u_{k+1})\psi_k\dx\,.
\end{equation}
We recognize two Riemann sums in the above expression, so that letting $n\to \infty$ gives
\[\begin{split}
\sum_{k=1}^{n-1}\int_{\Omega}(\psi_{k+1}-\psi_k)\AI u_k\dx
=\sum_{k=1}^{n-1}h\int_{\Omega}\frac{\psi_{k+1}-\psi_k}{h}\AI u_k\dx
\xrightarrow[n\to\, \infty]{}   \int_0^T\int_\Omega(\partial_t\psi)\,\AI (u)\dx\dt
\end{split}
\]
and
\[\begin{split}
\sum_{k=0}^{n-1}h\int_{\Omega}\n(u_{k+1})\psi_k\dx
\xrightarrow[n\to\, \infty]{}   \int_0^T \int_\Omega\n(u)\psi\dx\dt
\end{split}
\]
Therefore, taking limits as $n\to\infty$ in \eqref{1.3.Mild-WDS} gives:
\begin{equation}\label{1.6.Mild-WDS}
\lim_{n\to\infty}\int_{\Omega}\left(\psi_n \AI u_n -\psi_1\AI u_0\right)\dx =  \int_0^T\int_\Omega\AI (u)\partial_t\psi\dx\dt
-\int_0^T \int_\Omega\n(u)\psi\dx\dt\,.
\end{equation}
If we prove that the limit in the left hand side is zero, then $u$ would be a weak dual solution and the proof would be concluded.
Firstly, since by \eqref{psi.cpt} we have $\psi(t_1)= \psi(T/n)=0$ (at least for $n$ large), we easily get that
\[
\int_{\Omega} \psi_1\AI u_0 \dx=0\,.
\]
The same argument proves that
\[
\lim_{n\to\infty}\int_{\Omega}\psi_n \AI u_n\dx = \int_{\Omega}\psi(T,x) \AI u(T,x)\dx=0\,.
\]
We still have to check the way the initial data are taken. We know that mild solutions satisfy $u\in C([0,T]: \LL^1(\Omega))$. Since $\LL^1(\Omega)\subset\LL^1_{\p}(\Omega)$, we obtain the regularity that is needed in Definition \ref{Def.Very.Weak.Sol.Dual} of weak dual solutions to (CDP).\qed

\begin{cor}[Semigroup solutions with $u_0\in\LL^p$ are weak dual solutions]\label{Prop.VWSp-WDS}
Let $u$ be the unique semigroup (mild) solution corresponding to the initial datum $u_0\in \LL^p(\Omega)$ with $p\ge 1$\,, constructed in Theorem $\ref{thm.CP-JFA}$. Then $u$ is a weak dual solution in the sense of Definition $\ref{Def.Very.Weak.Sol.Dual}$ and is contained in the class $\mathcal{S}_p$\,.
\end{cor}

\noindent {\sl Proof.~}Since $u_0\in \LL^p(\Omega)$ with $p\ge 1$\,, then by Proposition \ref{Prop.Lp.Mild} we know that the unique semigroup (mild) solution $u(t)$ still belongs to $\LL^p(\Omega)$ for all $t>0$. Moreover, by Proposition \ref{Prop.Mild-WDS} we know that such semigroup solution is also a weak dual solution, which still belongs to $\LL^p(\Omega)$ for all $t\ge 0$, therefore is contained in the class $\mathcal{S}_p$.\qed

\section{Weighted $\LL^1$-estimates}\label{Sect.Weighted.L1}

As an interesting application of the weighted smoothing effects, we obtain the following $\LL^1$-weighted estimates, which will play an essential role in the existence proof, Theorem \ref{thm.L1weight.exist}. Though they are not contractivity statements, they will play the role of  weighted $\LL^1$-contractivity estimates when applied to ordered pairs of solutions.

In order to simplify the presentation, we first treat the case in which $\A$ has a first nonnegative eigenfunction $\Phi_1$; we recall that $\Phi_1\asymp \p$ on $\overline{\Omega}$, under the running assumption (K2).

\begin{prop}\label{thm.L1weight.contr}
Let $u\ge v$ be two ordered weak dual solutions to the Problem $(CDP)$ corresponding to the initial data $0\le u_0,v_0\in \LL^1_{\Phi_1}(\Omega)$\,. Then for all $t_1\ge t_0\ge 0$
\begin{equation}\label{L1weight.contr.estimates.0}
\int_{\Omega}\big[u(t_1,x)-v(t_1,x)\big]\Phi_1(x)\dx\le \int_{\Omega}\big[u(t_0,x)-v(t_0,x)\big]\Phi_1(x)\dx\,.
\end{equation}
Moreover, for all $0\le \tau_0\le \tau,t ~<+\infty$ such that either $t,\tau\le K_0$ or $\tau_0\ge K_0$\,, we have
\begin{equation}\label{L1weight.contr.estimates.1}\begin{split}
\int_{\Omega}\big[u(\tau,x)-v(\tau,x)\big]\Phi_1(x)\dx
&\le \int_{\Omega}\big[u(t,x)-v(t,x)\big]\Phi_1(x)\dx \\
&+ K_8[u(\tau_0)]\,\left|t-\tau\right|^{2s\vartheta_{i,\gamma}}\,\int_{\Omega} \big[u(\tau_0,x)-v(\tau_0,x)\big]\Phi_1\dx
\end{split}
\end{equation}
where $i=0$ if $t,\tau\le \| u(\tau_0)\|_{\LL^1_{\Phi_1}(\Omega)}^{\frac{2s}{d+\gamma}}$ and $i=1$ if $t,\tau\ge \| u(\tau_0)\|_{\LL^1_{\Phi_1}(\Omega)}^{\frac{2s}{d+\gamma}}$\,, and
\begin{equation}\label{L1weight.contr.estimates.1.a}
 K_8[u(\tau_0)]:=  \frac{\lambda_1
(\ka \vee 1)\,\n(K_7)}{2s\vartheta_{i,\gamma}}\|u(\tau_0)\|_{\LL^1_{\Phi_1}(\Omega)}^{2s(m_i-1) \vartheta_{i,\gamma}}
:=K_9\|u(\tau_0)\|_{\LL^1_{\Phi_1}(\Omega)}^{2s(m_i-1) \vartheta_{i,\gamma}}
\end{equation}
and $K_7>0$ is given in Theorem $\ref{thm.Upper.1.PME}$, $\ka$ in Lemma \eqref{Lem.N1.F}, and $K_0$ in Theorem $\ref{thm.Upper.PME}$\,.
\end{prop}

\noindent\textbf{Remark. }Formula \eqref{L1weight.contr.estimates.0} is the only part of the Theorem needed in the proof of the smoothing effects, Theorem \ref{thm.Upper.PME.Smoothing.F}, and its use implies no circularity.

\medskip

\noindent {\sl Proof of Proposition \ref{thm.L1weight.contr}. }We split it into several steps.

\noindent$\bullet~$\textsc{Step 1. }\textit{Monotonicity. }We begin by applying to $u$ and $v$ separately the definition of very weak solution \ref{Def.Very.Weak.Sol.Dual} in the form given in formula \eqref{step.1.thm.repr} of Step 1 of the proof of Proposition \ref{prop.point.est}, with the admissible test function $\psi=\Phi_1$ (recall that $\AI \Phi_1=\lambda_1^{-1}\Phi_1\ge 0$); therefore we get for any $t,t_0\ge 0$
\begin{equation}\label{L1weight.contr.Step.1.1}
\begin{split}
\int_{\Omega}\big[u(t_0,\cdot)-v(t_0,\cdot)\big]\Phi_1\dx- \int_{\Omega}\big[u(t,\cdot)-v(t,\cdot)\big]\Phi_1\dx
&=\lambda_1\int_{t_0}^{t_1}\int_{\Omega} \big(\n(u)-\n(v)\big)\Phi_1\dx\rd\tau\\
\end{split}
\end{equation}
From this equality, the fact that $F$ is increasing, and the fact that solutions are ordered, namely $u-v\ge 0$, it immediately follows the monotonicity property
\begin{equation}\label{L1weight.contr.Step.1.2}
\int_{\Omega}\big[u(t,x)-v(t,x)\big]\Phi_1(x)\dx\le\int_{\Omega}\big[u(\tau,x)-v(\tau,x)\big]\Phi_1(x)\dx\,,
\qquad\mbox{for all }0\le \tau\le t\,.
\end{equation}
The above inequality with $v=0$, implies the monotonicity of the $\LL^1_{\Phi_1}$ norm: $\|u(t)\|_{\LL^1_{\Phi_1}(\Omega)}\le \|u(\tau)\|_{\LL^1_{\Phi_1}(\Omega)}$\,, for all $0\le \tau\le t$.

\noindent$\bullet~$\textsc{Step 2. }We are going to prove the following inequality, valid for all $0\le v \le u$\,,
\begin{equation}\label{L1weight.contr.Step.2.1}
\n(u)-\n(v) \le K'_7\,(u-v)\,
\frac{\|u(\tau_0)\|_{\LL^1_{\Phi_1}(\Omega)}^{2s(m_i-1) \vartheta_{i,\gamma}}}{t^{(d+\gamma)(m_i-1)\vartheta_{i,\gamma}}}
\qquad\mbox{with}\qquad K'_7=(\ka \vee 1)\,\n(K_7)\,.
\end{equation}
with $i=1$ if $t\ge \| u(\tau_0)\|_{\LL^1_{\Phi_1}(\Omega)}^{\frac{2s}{d+\gamma}}$ and $i=0$ if $t\le \|u(\tau_0)\|_{\LL^1_{\Phi_1}(\Omega)}^{\frac{2s}{d+\gamma}}$.

Since $\n$ is convex, therefore $ \n(u)\le u \n'(u) \le \n(2u)$ for all $u\ge 0$\,, moreover,  Lemma \ref{Lem.N1.F} implies that for any $r_0\ge 0$\,, we have  $\n(x)\le \ka\vee 1\, \n(2r_0) (x/(2r_0))^{m_i}$ with $i=0$ if $x\in [0,2r_0]$ and $i=1$ if $x\ge 2r_0$\,. Therefore we have
\[
F'(U)\le \frac{F(2U)}{U} \le (\ka \vee 1)\, \frac{\n(r_0)}{r_0^{m_i-1}} U^{m_i-1}\,,
    \qquad\mbox{with $i=0$ if $U\in [0,r_0]$ and $i=1$ if $U\ge r_0$}\,.
\]
As a consequence, we obtain that for any $0\le v\le u\le U$:
\begin{equation}\label{L1weight.contr.Step.1.2a}
\n(u)-\n(v) \le \n'(u)(u-v)\le \n'(U)(u-v) \le (\ka \vee 1)\, \frac{\n(r_0)}{r_0^{m_i-1}} U^{m_i-1} (u-v)
\end{equation}
since $\n$ is convex, therefore the function $F'$ is non-decreasing.\\
Next, we recall the weighted smoothing effect \eqref{thm.Upper.PME.Smoothing.F}, namely for all $0\le \tau_0\le t$ we have
\begin{equation}\label{L1weight.contr.Step.1.2b}
\|u(t)\|_{\LL^\infty(\Omega)}
\le K_7\, \frac{\|u(\tau_0)\|_{\LL^1_{\Phi_1}(\Omega)}^{2s \vartheta_{j,\gamma}}}{t^{(d+\gamma)\vartheta_{j,\gamma}}}:=U_j\,,
\qquad\mbox{for all }t\ge 0\,,
\end{equation}
with $j=1$ if $t\ge \| u(\tau_0)\|_{\LL^1_{\Phi_1}(\Omega)}^{\frac{2s}{d+\gamma}}$ and $j=0$
if $t\le \| u(\tau_0)\|_{\LL^1_{\Phi_1}(\Omega)}^{\frac{2s}{d+\gamma}}$\,, which is equivalent to $j=1$ if $U_j\le K_7$ and $j=0$ if $U_j\ge K_7$. As a consequence, we can choose $r_0=K_7$ so that we can let $i=j$ and we can join \eqref{L1weight.contr.Step.1.2a} with \eqref{L1weight.contr.Step.1.2b} to obtain
\begin{equation}\label{L1weight.contr.Step.1.2c}
\n(u)-\n(v) \le (\ka \vee 1)\, \frac{\n(K_7)}{K_7^{m_i-1}} U_i^{m_i-1} (u-v)
=(\ka \vee 1)\, \frac{\n(K_7)}{K_7^{m_i-1}} K_7^{m_i-1}\,
\frac{\|u(\tau_0)\|_{\LL^1_{\Phi_1}(\Omega)}^{2s(m_i-1) \vartheta_{i,\gamma}}}{t^{(d+\gamma)(m_i-1)\vartheta_{i,\gamma}}}(u-v)
\end{equation}
with $j=1$ if $t\ge \| u(\tau_0)\|_{\LL^1_{\Phi_1}(\Omega)}^{\frac{2s}{d+\gamma}}$ and $j=0$ if $t\le \|u(\tau_0)\|_{\LL^1_{\Phi_1}(\Omega)}^{\frac{2s}{d+\gamma}}$.
This concludes the proof of \eqref{L1weight.contr.Step.2.1}.

\noindent$\bullet~$\textsc{Step 3. }The final step consists in applying the inequality of Step 2 to estimate the right hand side of formula \eqref{L1weight.contr.Step.1.1}. Let $0\le \tau_0\le t_0\le t$\,,
such that either $t_0,t_1\le \| u(\tau_0)\|_{\LL^1_{\Phi_1}(\Omega)}^{\frac{2s}{d+\gamma}}$ (hence $i=0$) or $t_1,t_0\ge \| u(\tau_0)\|_{\LL^1_{\Phi_1}(\Omega)}^{\frac{2s}{d+\gamma}}$ (hence $i=1$)\,. As a consequence, in what follows the index $i$ will not change. Then we can estimate
\begin{equation}\label{L1weight.contr.Step.1.3}\begin{split}
\int_{t_0}^{t}\int_{\Omega} & \big(\n(u)-\n(v)\big)\Phi_1\dx\rd\tau
\le K'_7 \int_{t_0}^{t}
\frac{\|u(\tau_0)\|_{\LL^1_{\Phi_1}(\Omega)}^{2s(m_i-1) \vartheta_{i,\gamma}}}{\tau^{(d+\gamma)(m_i-1)\vartheta_{i,\gamma}}}
    \int_{\Omega} \big(u(\tau,x)-v(\tau,x)\big)\Phi_1(x)\dx\rd\tau\\
&\le K'_7 \|u(\tau_0)\|_{\LL^1_{\Phi_1}(\Omega)}^{2s(m_i-1) \vartheta_{i,\gamma}} \frac{(t-t_0)^{2s\vartheta_{i,\gamma}}}{2s\vartheta_{i,\gamma}}\int_{\Omega} \big(u(\tau_0,x)-v(\tau_0,x)\big)\Phi_1(x)\dx
\end{split}\end{equation}
where in the last step we have used the monotonicity inequality \eqref{L1weight.contr.Step.1.2} of Step 1\,, together with
\[
\int_{t_0}^{t}\frac{1}{\tau^{(d+\gamma)(m_i-1)\vartheta_{i,\gamma}}} \rd\tau = \frac{1}{2s\vartheta_{i,\gamma}}(t^{2s\vartheta_{i,\gamma}}-t_0^{2s\vartheta_{i,\gamma}})ç
\le \frac{(t-t_0)^{2s\vartheta_{i,\gamma}}}{2s\vartheta_{i,\gamma}}
\]
the last step is valid since ${2s\vartheta_{i,\gamma}}\le 1$\,. Plugging inequality \eqref{L1weight.contr.Step.1.3} into \eqref{L1weight.contr.Step.1.1} gives for all $t,t_0\ge \tau$\,:
\begin{equation}\label{L1weight.contr.Step.1.4}
\begin{split}
&\left| \int_{\Omega}\big[u(t_0,x)-v(t_0,x)\big]\Phi_1(x)\dx - \int_{\Omega}\big[u(t,x)-v(t,x)\big]\Phi_1(x)\dx\right|\\
&\qquad\qquad\le \lambda_1
K'_7 \|u(\tau_0)\|_{\LL^1_{\Phi_1}(\Omega)}^{2s(m_i-1) \vartheta_{i,\gamma}} \frac{|t-t_0|^{2s\vartheta_{i,\gamma}}}{2s\vartheta_{i,\gamma}}\int_{\Omega} \big(u(\tau_0,x)-v(\tau_0,x)\big)\Phi_1(x)\dx\\
\end{split}
\end{equation}
which implies \eqref{L1weight.contr.estimates.1}. The constant $K_8[u(\tau_0)]=\lambda_1
K'_7 \|u(\tau_0)\|_{\LL^1_{\Phi_1}(\Omega)}^{2s(m_i-1) \vartheta_{i,\gamma}} /(2s\vartheta_{i,\gamma})$\,, has the form given in \eqref{L1weight.contr.estimates.1.a}.\qed

\subsection{More general operators}

In the case in which $\A$ does not have a first eigenfunction $\Phi_1$ the results of Proposition \ref{thm.L1weight.contr} continue to hold for different weights. Indeed we observe that taking any nonnegative function $\psi\in \LL^\infty(\Omega)$\,, thanks to assumption (K2) it is easy to show that $\AI\psi\ge 0$ and
\begin{equation}\label{weight.1}
\AI\psi(x)\asymp \p(x)\qquad\mbox{for a.e. }x\in \Omega\,.
\end{equation}
This will imply the monotonicity of some $\LL^1$-weighted norm. Let's put $\Psi_1= \AI\p$\,, which bears some similarity with the formula $\Phi_1=\lambda_1^{-1}\AI\Phi_1$.
\begin{prop}\label{thm.L1weight.contr.psi}
Let $u\ge v$ be two ordered weak dual solutions to the Problem $(CDP)$ corresponding to the initial data $0\le u_0,v_0\in \LL^1_{\p}(\Omega)$\,. Then for all $0\le \psi\in \LL^\infty(\Omega)$
\begin{equation}\label{monotonicity.psi}
\int_{\Omega}\big[u(t,x)-v(t,x)\big]\AI\psi(x)\dx\le\int_{\Omega}\big[u(\tau,x)-v(\tau,x)\big]\AI\psi(x)\dx\,,
\qquad\mbox{for all }0\le \tau\le t\,.
\end{equation}
As a consequence, there exists a constant $C_{\Omega,\gamma}>0$ such that
\begin{equation}\label{monotonicity.dist}
\int_{\Omega}\big[u(t,x)-v(t,x)\big]\,\p(x)\dx\le C_{\Omega,\gamma}\int_{\Omega}\big[u(\tau,x)-v(\tau,x)\big]\,\p(x)\dx\,,
\qquad\mbox{for all }0\le \tau\le t\,.
\end{equation}
Moreover, for all $0\le \tau_0\le \tau,t ~<+\infty$ such that either $t,\tau\le K_0$ or $\tau_0\ge K_0$\,, we have
\begin{equation}\label{L1weight.contr.estimates.1.dist}\begin{split}
\int_{\Omega}\big[u(\tau,x)-v(\tau,x)\big] \Psi_1(x)\dx
&\le \int_{\Omega}\big[u(t,x)-v(t,x)\big] \Psi_1(x)\dx \\
&+ K_8[u(\tau_0)]\,\left|t-\tau\right|^{2s\vartheta_{i,\gamma}}\,\int_{\Omega} \big[u(\tau_0,x)-v(\tau_0,x)\big]\p(x)\dx
\end{split}
\end{equation}
where $i=0$ if $t,\tau\le \| u(\tau_0)\|_{\LL^1_{\p}(\Omega)}^{\frac{2s}{d+\gamma}}$ and $i=1$ if $t,\tau\ge \| u(\tau_0)\|_{\LL^1_{\p}(\Omega)}^{\frac{2s}{d+\gamma}}$\,, and
\begin{equation}\label{L1weight.contr.estimates.1.ab}
 K_8[u(\tau_0)]:=  \frac{\lambda_1
(\ka \vee 1)\,\n(K_7)}{2s\vartheta_{i,\gamma}}\|u(\tau_0)\|_{\LL^1_{\p}(\Omega)}^{2s(m_i-1) \vartheta_{i,\gamma}}
:=K_9\|u(\tau_0)\|_{\LL^1_{\p}(\Omega)}^{2s(m_i-1) \vartheta_{i,\gamma}}
\end{equation}
and $K_7>0$ is given in Theorem $\ref{thm.Upper.1.PME}$, $\ka$ in Lemma \eqref{Lem.N1.F}, and $K_0$ in Theorem $\ref{thm.Upper.PME}$\,.
\end{prop}

\noindent {\sl Proof.~}We split several steps.

\noindent$\bullet~$\textsc{Step 1. }\textit{Monotonicity. }We begin by applying to $u$ and $v$ separately the definition of very weak solution \ref{Def.Very.Weak.Sol.Dual} in the form given in formula \eqref{step.1.thm.repr} of Step 1 of the proof of Proposition \ref{prop.point.est}, with the admissible test function $0\le \psi \in \LL^\infty(\Omega)$ (recall that $\AI \psi\ge 0$); therefore we get for any $t\ge t_0\ge 0$
\begin{equation}\label{L1weight.contr.Step.1.1.psi}
\begin{split}
\int_{\Omega}\big[u(t_0,x)-v(t_0,x)\big]\AI\psi(x)\dx &- \int_{\Omega}\big[u(t,x)-v(t,x)\big]\AI\psi(x)\dx\\
&=\int_{t_0}^{t}\int_{\Omega} \big[\n(u(\tau,x))-\n(v(\tau,x))\big]\psi(x)\dx\rd\tau\ge 0\,.\\
\end{split}
\end{equation}
From this equality, the fact that $F$ is increasing, and the fact that solutions are ordered, namely $u-v\ge 0$, it immediately follows the monotonicity property \eqref{monotonicity.psi}. Notice that letting $v=0$, gives the monotonicity of the $\LL^1_{\AI\psi}$ norm of $u(t,\cdot)$\,: $\|u(t)\|_{\LL^1_{\AI\psi}(\Omega)}\le \|u(\tau)\|_{\LL^1_{\AI\psi}(\Omega)}$\,, for all $0\le \tau\le t$.

\noindent Finally, as a consequence of  \eqref{monotonicity.psi} and \eqref{weight.1} we get \eqref{monotonicity.dist}\,.

\noindent$\bullet~$\textsc{Step 2 and 3. }We can repeat Step 2 and 3 of the proof of Proposition \eqref{thm.L1weight.contr}, to estimate the right-hand side of \eqref{L1weight.contr.Step.1.1.psi}, just by replacing $\Phi_1$ with $\psi$, we therefore get
\begin{equation*}
\begin{split}
\int_{t_0}^{t}\int_{\Omega} \big(\n(u)-\n(v)\big)\psi\dx\rd\tau
&\le K'_7 \|u(\tau_0)\|_{\LL^1_{\p}(\Omega)}^{2s(m_i-1) \vartheta_{i,\gamma}} \frac{(t-t_0)^{2s\vartheta_{i,\gamma}}}{2s\vartheta_{i,\gamma}}\int_{\Omega} \big(u(\tau_0,x)-v(\tau_0,x)\big)\psi(x)\dx
\end{split}\end{equation*}
which combined with \eqref{L1weight.contr.Step.1.1.psi}, gives \eqref{L1weight.contr.estimates.1.dist}, once we recall that $\AI\p\asymp \p$.\qed


\section{Proof of existence and uniqueness of weak dual solutions}\label{sect.exist.uniq.WDS}

\noindent{\bf Proof of Theorem \ref{thm.L1weight.exist} (Existence)} \\
\noindent$\bullet~$\textsc{Step 1.~}\textit{Construction of limit solutions via monotone sequences of bounded mild solutions. }We consider general data $u_0\ge 0$\,, and construct a limit solution by approximation from below with $\LL^\infty$-mild solutions. We consider a monotone non-decreasing sequence $0\le u_{0,n}\le u_{0,n+1}\le u_0$\,, with $u_{0,n}\in \LL^\infty(\Omega)$\,, monotonically converging from below to $u_0\in \LL^1_{\p}(\Omega)$ in the topology of $\LL^1_{\p}(\Omega)$. By Theorem \ref{thm.CP-JFA}, we know that to every $u_{0,n}$ corresponds a unique mild solution $u_n(t,x)\in \LL^\infty(\Omega)$ for all $t>0$, by Proposition \ref{Prop.Lp.Mild}. Next, by Proposition \ref{Prop.Mild-WDS} we know that mild solutions are weak dual solutions in the sense of Definition \ref{Def.Very.Weak.Sol.Dual}\,.  Moreover, since comparison holds for mild solutions (cf. Theorem \ref{thm.CP-JFA}), then the sequence $u_n(t,x)$ is ordered, namely $u_n(t,x)\le u_{n+1}(t,x)$ for all $x\in \Omega$ and $t>0$\,. We have that $u_n\in \mathcal{S}_\infty$\,, therefore all the upper estimates of Theorems \ref{thm.Upper.PME} and \ref{thm.Upper.PME.Smoothing.F} together with their corollaries, hold true for $u_n$\,. In particular, since the absolute upper bounds of Theorem \eqref{thm.Upper.PME} are independent on the initial datum, we can guarantee that for any fixed $\tau>0$, there exists the monotone limit $u(t,x)=\lim\limits_{n\to \infty} u_n(t,x)$ in $\LL^\infty((\tau,\infty)\times\Omega)$, and that such limit satisfies the same upper estimates of $u_n$, by the lower semicontinuity of the $\LL^\infty$ norm.

\noindent$\bullet~$\textsc{Step 2.~}\textit{Next, we prove that  $u\in C^0([0,\infty)\,:\,\LL^1_{\p}(\Omega))$. }On one hand, by monotonicity estimates \eqref{CP.Monotonicity.3}, namely that $t\mapsto t^{\frac{1}{m_0-1}}\,u(t,x)$ is nondecreasing in $t>0$ for a.e. $x\in \Omega$\,, which implies that for all $0\le t_0\le t_1$
\begin{equation}\label{cont.0}
u(t_1,x)- u(t_0,x)\ge \left(\frac{t_0}{t_1}\right)^{\frac{1}{m_0-1}}u(t_0,x)-u(t_0,x)
=-\left[1- \left(\frac{t_0}{t_1}\right)^{\frac{1}{m_0-1}}\right]u(t_0,x)
\end{equation}
which implies that the negative part $(u(t_1,\cdot)- u(t_0,\cdot))_-$ satisfies
\begin{equation}\label{cont.1}
0\le (u(t_1,x)- u(t_0,x))_- \le \left[1- \left(\frac{t_0}{t_1}\right)^{\frac{1}{m_0-1}}\right]u(t_0,x)\qquad\mbox{for all }0\le t_0\le t_1\,.
\end{equation}
Integrating the above inequality with the weight $\Psi_1=\AI\p$\,, gives for all $0\le t_0\le t_1$
\begin{equation}\label{cont.2}
0\le\int_\Omega (u(t_1,x)- u(t_0,x))_-\, \Psi_1(x)\dx\le \left[1- \left(\frac{t_0}{t_1}\right)^{\frac{1}{m_0-1}}\right]\int_\Omega u(t_0,x)\Psi_1(x)\dx
\end{equation}

On the other hand, using the estimates \eqref{monotonicity.dist} for $u_n\le u$, we have that for all $\tau\ge 0$
\begin{equation}\label{Step.1.1.thm.L1weight.exist}\begin{split}
0\le \int_{\Omega}\big[u(\tau,x)-u_n(\tau,x)\big]\p(x)\dx
&\le C_{\Omega,\gamma}\int_{\Omega}\big[u_0(x)-u_{0,n}(x)\big]\p(x)\dx \\
\end{split}
\end{equation}
hence $u_n(\tau)\to u(\tau)$ for all $\tau\ge 0$ as $n\to \infty$ in the strong $\LL^1_{\p}$ topology and we can pass to the limit in inequality \eqref{L1weight.contr.estimates.1.dist} (with $v=0$) to obtain for all $t_1\ge t_0\ge 0$

\begin{equation}\label{cont.3}\begin{split}
\left|\int_{\Omega}u(t_1,x)\Psi_1(x)\dx
-\int_{\Omega}u(t_0,x)\Psi_1(x)\dx \right|
\le K_8[u_0]\,\left|t_1-t_0\right|^{2s\vartheta_{i,\gamma}}\,\int_{\Omega} u(t_0,x)\Psi_1(x)\dx
\end{split}
\end{equation}
Recalling that $|f|=f+2f_-$\,, with $f_-=\max\{0,-f\}$, and letting $f=u(t_1,x)- u(t_0,x)$ we obtain, for all $t_1\ge t_0\ge 0$\,:
\begin{equation}\label{cont.4}\begin{split}
\int_{\Omega}|u(t_1,x)- u(t_0,x)| & \Psi_1(x)\dx\\
    &=\int_{\Omega}(u(t_1,x)- u(t_0,x))\Psi_1(x)\dx
        +2\int_\Omega (u(t_1,x)- u(t_0,x))_-\, \Psi_1(x)\dx\\
    &\le \left[K_8[u_0]\,\left|t_1-t_0\right|^{2s\vartheta_{i,\gamma}} +2\left(1- \left(\frac{t_0}{t_1}\right)^{\frac{1}{m_0-1}}\right) \right]\int_\Omega u(t_0,x)\Psi_1(x)\dx\\
\end{split}
\end{equation}
where in the last step we have used inequalities \eqref{cont.2} and \eqref{cont.3}.

It remains to prove the continuity at $t=0$.
\begin{equation}\begin{split}
\int_{\Omega}|u(t,x)- u_0(x)|  \Psi_1(x)\dx
&\le \int_{\Omega}|u(t,x)- u_n(t,x)|  \Psi_1(x)\dx
+   \int_{\Omega}|u_n(t,x)-u_{0,n}(x)|  \Psi_1(x)\dx\\
&+   \int_{\Omega}|u_{0,n}(x)-u_0|  \Psi_1(x)\dx= (I)+(II)+(III)
\end{split}\end{equation}
Let us fix $\varepsilon >0$, then $(III)\le \varepsilon/3$ if $n$ is large, by construction.
Let us estimate the remaining two terms as follows. For the first term we use inequality \eqref{monotonicity.psi}, since by construction $u\ge u_n$, so that
\[
\int_{\Omega}|u(t,x)- u_n(t,x)|  \Psi_1(x)\dx
\le \int_{\Omega}|u_0(x)- u_{0,n}(x)|  \Psi_1(x)\dx= (III)
\]
As for the second term we have
\[
\int_{\Omega}|u_n(t,x)-u_{0,n}(x)|  \Psi_1(x)\dx
\le \| \Psi_1\|_{\LL^\infty(\Omega)} \int_{\Omega}|u_n(t,x)-u_{0,n}(x)| \dx\,.
\]
and we recall that $n$ is fixed,  $u_n\in C^0([0,\infty)\,:\,\LL^1(\Omega))$, so that the above expression goes to zero as $t\to 0$.

The above estimates show that $u\in C^0([0,\infty)\,:\,\LL^1_{\Psi_1}(\Omega))$. Then, since $\Psi_1=\AI\p\asymp\p$ then the weighted norms with $\Psi_1$ and $\p$ are equivalent and this is enough to conclude that $u\in C^0([0,\infty)\,:\,\LL^1_{\p}(\Omega))$.

As a final remark, we have showed that the solution constructed by approximation has the required properties to be in the class $\mathcal{S}_\infty$, provided we prove that it is a weak dual solution in the sense of Definition \ref{Def.Very.Weak.Sol.Dual}; the latter will be proved in Step 3.

\noindent$\bullet~$\textsc{Step 3. }\textit{The limit solution is a weak dual solution. }We have to check that for all $\psi$ such that $\psi/\p\in C^1_c((0,+\infty): \LL^\infty(\Omega))$ the following identity holds true:
\begin{equation}
\displaystyle \int_0^\infty\int_{\Omega}\AI (u) \,\dfrac{\partial \psi}{\partial t}\,\dx\dt
-\int_0^\infty\int_{\Omega} \n(u)\,\psi\,\dx \dt=0.
\end{equation}
Let us fix an admissible test function $\psi$\,. Proposition \ref{Prop.Mild-WDS} shows that mild solutions are weak dual solutions in the sense of Definition \ref{Def.Very.Weak.Sol.Dual}\,,  so that
\begin{equation}
\displaystyle \int_0^\infty\int_{\Omega}\AI (u_n) \,\dfrac{\partial \psi}{\partial t}\,\dx\dt
=\int_0^\infty\int_{\Omega}\n(u_n)\,\psi\,\dx \dt\,.
\end{equation}
The proof is concluded once we show that
\begin{equation}\label{Step.2.1.thm.L1weight.exist}
\int_0^\infty\int_{\Omega}\AI (u_n) \,\dfrac{\partial \psi}{\partial t}\,\dx\dt
\xrightarrow[n\to\, \infty]{}
\int_0^\infty\int_{\Omega}\AI (u) \,\dfrac{\partial \psi}{\partial t}\,\dx\dt\,,
\end{equation}
and also that
\begin{equation}\label{Step.2.2.thm.L1weight.exist}
\int_0^\infty\int_{\Omega} \n(u_n)\,\psi\,\dx \dt
\xrightarrow[n\to\, \infty]{}
\int_0^\infty\int_{\Omega} \n(u)\,\psi\,\dx \dt\,.
\end{equation}
\noindent$\circ~$\textit{Proof of \eqref{Step.2.1.thm.L1weight.exist}. }Recall that $\psi/\p\in C^1_c((0,+\infty): \LL^\infty(\Omega))$ and say that the time-support is contained in $[t_1,t_2]$\,, and recall that $f\ge 0$ implies $\AI f \ge 0$\,, so that $u\ge u_n$ implies $\AI(u-u_n)\ge 0$ and
\begin{equation*}\begin{split}
& \left|\int_0^\infty\int_{\Omega}\AI (u-u_n)(t,x) \,\partial_t\psi\,\dx\dt\right|
\le \int_{t_1}^{t_2} \left\|\frac{\partial_t\psi}{\p}\right\|_{\LL^\infty(\Omega)}\int_{\Omega}\p(x)\AI (u-u_n)(t,x) \,\dx\dt\\
&= \int_{t_1}^{t_2} \left\|\frac{\partial_t\psi}{\p}\right\|_{\LL^\infty(\Omega)}\int_{\Omega}\AI\p(x) (u-u_n)(t,x) \,\dx\dt
\le C_{\Omega,\gamma}\int_{t_1}^{t_2} \left\|\frac{\partial_t\psi}{\p}\right\|_{\LL^\infty(\Omega)}\int_{\Omega} (u-u_n)(t,x)\p(x) \,\dx\dt\\
& \le C_{\Omega,\gamma}\int_{\Omega} (u_0-u_{0,n})\p\,\dx\, \int_{t_1}^{t_2} \left\|\frac{\partial_t\psi}{\p}\right\|_{\LL^\infty(\Omega)}\dt \xrightarrow[n\to\, \infty]{} 0 \\
\end{split}
\end{equation*}
where we have used the fact that $\AI$ is symmetric and in the last step we have used inequality \eqref{Step.1.1.thm.L1weight.exist}\,.

\noindent$\circ~$\textit{Proof of \eqref{Step.2.2.thm.L1weight.exist}. }Recall that $\psi/\p\in C^1_c((0,+\infty): \LL^\infty(\Omega))$\,, and say that the time-support is contained in $[t_1,t_2]\subset (0,\infty)$\,, so that $\left\|\psi/\p\right\|_{\LL^\infty([t_1,t_2]\times\Omega)}\le \ka_1$\,, and
\begin{equation*}\begin{split}
\left|\int_0^\infty\int_{\Omega} (\n(u)-\n(u_n))\,\psi\,\dx \dt\right|
&\le \left\|\frac{\psi}{\p}\right\|_{\LL^\infty([t_1,t_2]\times\Omega)}\int_{t_1}^{t_2} \int_{\Omega}(\n(u)-\n(u_n))(t,x)\,\p(x) \,\dx\dt\\
&\le \ka_1 \,\ka_2[u_0](t_1-t_2)^{2s\vartheta_{i,\gamma}}\,\int_{\Omega}(u_0-u_{0,n})\,\p \,\dx
\xrightarrow[n\to\, \infty]{} 0 \\
\end{split}
\end{equation*}
where in the last step we have used inequality \eqref{L1weight.contr.Step.1.3}, namely for all $0\le t_1\le t_2$
\begin{equation}\label{thm.esist.01}\begin{split}
\int_{t_1}^{t_2}\int_{\Omega} \big(\n(u)-\n(v)\big)\p\dx\rd\tau
&\le K'_7 \|u_0\|_{\LL^1_{\p}(\Omega)}^{2s(m_i-1) \vartheta_{i,\gamma}} \frac{(t_1-t_2)^{2s\vartheta_{i,\gamma}}}{2s\vartheta_{i,\gamma}}\int_{\Omega} \big(u_0(x)-u_{0,n}(x)\big)\p(x)\dx\\
&:=\ka_2[u_0]\, (t_1-t_2)^{2s\vartheta_{i,\gamma}}\,\int_{\Omega} \big(u_0(x)-u_{0,n}(x)\big)\p(x)\dx\mbox{\,.\qed}
\end{split}\end{equation}

\medskip

\noindent{\bf Proof of Theorem \ref{thm.L1weight.uniq} (Uniqueness)} \\
\noindent {\sl Proof.~} We keep the notations of the proof of Theorem \ref{thm.L1weight.exist}.
Assume that there exist another monotone non-decreasing sequence $0\le v_{0,k}\le v_{0,k+1}\le u_0$\,, with $v_{0,k}\in \LL^\infty(\Omega)$\,, monotonically converging from below to $u_0\in \LL^1_{\p}(\Omega)$ in the topology of $\LL^1_{\p}(\Omega)$. By the same considerations as in the proof of Theorem \ref{thm.L1weight.exist} we can show that there exists a solution $v(t,x)\in C^0([0,\infty)\,:\,\LL^1_{\p}(\Omega))$. We want to show that $u=v$, where $u$ is the solution constructed in the same way from the sequence $u_{0,n}$. We will prove equality by proving that $v\le u$ and then that $u\le v$. To prove that $v\le u$ we use the estimates
\begin{equation}\label{final.001}
\int_{\RR^d}\big[v_k(t,x)-u_n(t,x)\big]_+\dx\le \int_{\RR^d}\big[v_k(0,x)-u_n(0,x)\big]_+\dx
\end{equation}
which hold for all $u_n(t,\cdot)$ and  $v_k(t,\cdot)$, since both are mild solutions constructed in Theorem \cite{CP-JFA}.\\
Letting $n\to \infty$ we get that
\[
\lim_{n\to\infty}\int_{\RR^d}\big[v_k(t,x)-u_n(t,x)\big]_+\dx
\le \lim_{n\to\infty}\int_{\RR^d}\big[v_k(0,x)-u_n(0,x)\big]_+\dx
=\int_{\RR^d}\big[v_k(0,x)-u_0(x)\big]_+\dx=0
\]
since $v_k(0,x)\le u_0$ by construction. Therefore also $v_k(t,x)\le u(t,x)$ for $t>0$, so that in the limit $k\to \infty$ we obtain $v(t,x)\le u(t,x)$\,. The inequality $u\le v$ can be obtained simply by switching the roles of $u_n$ and $v_k$\,. The validity of estimates of Proposition $\ref{thm.L1weight.contr.psi}$ is guaranteed by the above limiting process. The comparison holds by taking the limits in inequality \eqref{final.001}\,.\qed

\section{Appendix}

\subsection{Technical Proofs}
Put here the proof of the Pointwise estimates, and of the numerical lemmas for F and its Legendre transform.

\subsubsection{Proof of Proposition \ref{prop.point.est}}\label{app.proof.5.1}

\noindent\textsc{Proof of Proposition \ref{prop.point.est}. }We adapt the scheme of the proof of Proposition 4.2 of \cite{BV-PPR1}, in which we have treated the case $\n(u)=u^m$. The proof consists of several steps.

\noindent \textsc{Step 1. }\textit{$\LL^1$-weighted estimates. } We will use the definition \eqref{Def.Very.Weak.Sol.Dual} of weak dual solution, with a test function of the form $\psi(t,x)=\psi_1(t)\psi_2(x)$, where $\psi_1(t)\in  C^1_c((0,+\infty))$ and $\psi_2/\p\in \LL^\infty(\Omega)$. It follows that $u\in C((0,\infty): \LL^1_{\p}(\Omega))$, $\n(u) \in \LL^1\left((0,\infty):\LL^1_{\p}(\Omega)\right)$ satisfies the identity
\begin{equation}\label{step.1.thm.repr.0}\begin{split}
\int_0^{+\infty}\psi'_1(\tau)\int_{\Omega}u(\tau,x)  \AI\psi_2(x)\dx\rd\tau
 &= \int_0^{+\infty}\psi_1(\tau)\int_{\Omega}\n(u(\tau,x))\,\psi_2(x)\dx \rd\tau\,,
\end{split}
\end{equation}
where in the left-hand side we have used the symmetry of the operator $\AI$. Notice that the space integral on the left-hand side of the formula is bounded, because of the following argument: we write $\psi_2= v\p$ with $v(x)$ bounded, and recall that $u\ge 0$. Then,
\[
\int_\Omega \AI u(t)\, \psi_2\dx \le \left\|v\right\|_{\LL^\infty}\int_\Omega \p  \AI u(t)\dx
=\left\|v\right\|_{\LL^\infty}\int_\Omega u(t)\AI \p\dx\le C_{\Omega,\gamma}\left\|v\right\|_{\LL^\infty} \int_\Omega u(t)\p\dx\,.
\]
where in the last step we have used the fact that $\AI\p\asymp\p$\,, which follows by assumption (K2).
We now want to pass to the limit in \eqref{step.1.thm.repr.0} and prove that for all $0\le t_0\le t_1$ and for all $\psi_2(x)$\,, with $\psi_2: \overline{\Omega}\to \RR$ measurable and $\|\psi_2/\p\|_{\LL^\infty(\Omega)}<+\infty$\,, we have
\begin{equation}\label{step.1.thm.repr}\begin{split}
\int_{\Omega}u(t_0,x)\AI \psi_2(x)\dx - \int_{\Omega}u({t_1},x)\AI \psi_2(x)\dx=\int_{t_0}^{t_1}\int_{\Omega}\n(u(\tau,x))\psi_2(x)\dx \rd\tau\,.
\end{split}
\end{equation}
This is rather standard: we only need to take $\psi_1(\tau)=\chi_{[t_0,t_1]}(\tau)$ as test function in formula \eqref{step.1.thm.repr.0}\,, so that $\psi_1'(\tau)=\delta_{t_0}(\tau)-\delta_{t_1}(\tau)$; this can be jusfified by considering a smooth approximation $\psi_{1,n}\in C_c^{\infty}(0,+\infty)$ so that $\psi_{1,n}\to \chi_{[t_0,t_1]}(\tau)$ in $\LL^\infty(0,+\infty)$\,, and so that $\psi'_{1,n} \to \delta_{t_0}(\tau)-\delta_{t_1}(\tau)$ in the sense of Radon measures with compact support. Clearly, these approximations are admissible test functions such that $\psi_{n}/\p\in C^1_c((0,+\infty): \LL^\infty(\Omega))$\,. Under the above assumptions, it is clear that
\[
\int_0^{+\infty}\psi'_{1,n}(\tau)\int_{\Omega}u(\tau,x)  \AI\psi_2(x)\dx\rd\tau \xrightarrow[n\to\, \infty]{}\int_{\Omega}u(t_0,x)\AI \psi_2(x)\dx -\int_{\Omega}u({t_1},x)\AI \psi_2(x)\dx\,,
\]
since $u\in C((0,\infty): \LL^1_{\p}(\Omega))$ implies that $\int_{\Omega}u({t_1},x)\AI \psi_2(x)\dx \in C^0(0,\infty)$\,, and since $\psi'_{1,n} \to \delta_{t_0}(\tau)-\delta_{t_1}(\tau)$  in the sense of Radon measures with compact support. On the other hand,
\[
\int_0^{+\infty}\psi_{1,n}(\tau)\int_{\Omega}\n(u(\tau,x))\,\psi_2(x)\dx \rd\tau \xrightarrow[n\to\, \infty]{}
\int_{t_0}^{t_1}\int_{\Omega}\n(u(\tau,x))\psi_2(x)\dx \rd\tau
\]
since $\n(u) \in \LL^1\left((0,\infty):\LL^1_{\p}(\Omega)\right)$ implies that $\int_{\Omega}\n(u(\tau,x))\psi_2(x)\dx \in \LL^1(0,\infty)$\,, and $\psi_{1,n}\to \chi_{[t_0,t_1]}(\tau)$ in $\LL^\infty(0,+\infty)$\,.

\noindent$\bullet~$\noindent\textsc{Step 2. }\textit{Proof of \eqref{thm.NLE.PME.estim.0}. }From \eqref{step.1.thm.repr} we  prove estimate \eqref{thm.NLE.PME.estim.0}, by first fixing $x_0\in \Omega$ and then taking a sequence of nonnegative test functions $\psi_{2,n}^{(x_0)}$ with $\psi_{2,n}^{(x_0)}: \overline{\Omega}\to \RR$ measurable and $\|\psi_2/\Phi_1\|_{\LL^\infty(\Omega)}<+\infty$\,, such that $\psi_{2,n}^{(x_0)}\to \delta_{x_0}$ as $n\to \infty$\,, in the sense of Radon measures. Therefore, $\AI \psi_{2,n}^{(x_0)}\to \K(\cdot , x_0)$ in $\LL^q(\Omega)$ for all $0<q<N/(N-2s)$, so that, taking $p=q/(q-1)>N/(2s)$\,, we have that
\[
\left|\int_{\Omega}u(\tau,x)  \AI\psi_{2,n}(x)\dx-\int_{\Omega}u(\tau,x) \K(x, x_0)\dx\right|
\le \|u(\tau)\|_{\LL^p(\Omega)}\|\AI\psi_{2,n}-\K(\cdot , x_0)\|_{\LL^q(\Omega)}\to 0
\]
as $n\to\infty$\,,  and we recall that  $u\in \mathcal{S}$\,, therefore $u(t)\in \LL^p(\Omega)$ for all $t> 0$\,, with $p>N/(2s)$\,. Since the right-hand side of \eqref{step.1.thm.repr} is non-negative, we have proved the first estimate of the Theorem, \eqref{thm.NLE.PME.estim.0}\,.

\noindent$\bullet~$\noindent\textsc{Step 3. }\textit{Monotonicity estimates} We will use  estimates \eqref{CP.Monotonicity} of \cite{CP-JFA}, namely that the function
\[
t\mapsto t^{\frac{1}{\mu_0}}\,\n(u(t,x))\qquad\mbox{is nondecreasing in $t>0$ for a.e. }x\in \Omega\,.
\quad\left(\mbox{Recall also that }\frac{1}{\mu_0}=\frac{m_0}{m_0-1}\right)
\]
to estimate the right-hand side of identity \eqref{step.1.thm.repr} from below and from above. More precisely, we will prove the following

\noindent {\bf Claim:} \textsl{For  almost every $0\le t_0\le t_1 \le t$\,, and for all $\psi_2$ as in Step 1, we have:}
\begin{equation}\label{step.2.1.thm.repr}
\begin{split}
\left(\frac{t_0}{t_1}\right)^{\frac{m_0}{m_0-1}}(t_1-t_0)\int_\Omega \n(u(t_0,x))\psi_2(x)\dx
&\le \int_{t_0}^{t_1}\int_\Omega \n(u(\tau,x)) \psi_2(x)\dx \rd\tau \\
&\le \frac{m_0-1}{t_0^{\frac{1}{m_0-1}}} t^{\frac{m_0}{m_0-1}}\,\int_\Omega \n(u(t,x))\psi_2(x)\dx\,.
\end{split}
\end{equation}
Consider a smooth sequence  $\psi_{1,n}\in C_c^{\infty}(0,+\infty)$, $0\le \psi_{1,n}\le 1$, such that $\psi_{1,n}\to \chi_{[t_0,t_1]}$ in $\LL^\infty(0,+\infty)$ and such that $\supp(\psi_{1,n})\subseteq [t_0-1/n, t_1+1/n]$ and $\psi_{1,n}\ge \chi_{[t_0,t_1]}$\,. There are two cases.

\medskip

\noindent\textit{Upper estimates. }Let $n$ be so big that $0\le t_0-1/n \le t_1+1/n \le t$\,, and recall that $u\ge 0$\,, so that
\[
\begin{split}
\int_{0}^{\infty}\psi_{1,n}(\tau) &\int_{\Omega}\n(u(\tau,x))\psi_2(x)\dx \rd\tau
\le  \int_{0}^{\infty}\psi_{1,n}(\tau) \left(\frac{t}{\tau}\right)^{\frac{m_0}{m_0-1}}\rd\tau \, \int_{\Omega}\n(u(t ,x))\psi_2(x)\dx \\
&\le \|\psi_{1,n}\|_{\LL^\infty(0,+\infty)} \int_{t_0-1/n}^{t_1+1/n}\left(\frac{t }{\tau}\right)^{\frac{m_0}{m_0-1}}\rd\tau \, \int_{\Omega}\n(u(t,x))\psi_2(x)\dx \\
&=(m_0-1)t^{\frac{m_0}{m_0-1}} \|\psi_{1,n}\|_{\LL^\infty(0,+\infty)} \left[\left(\frac{1}{t_0-\frac{1}{n}}\right)^{\frac{1}{m_0-1}}-\left(\frac{1}{t_1+\frac{1}{n}}\right)^{\frac{1}{m_0-1}}\right]\, \int_{\Omega}\n(u(t,x))\psi_2(x)\dx \\
&\le (m_0-1)\frac{\|\psi_{1,n}\|_{\infty} }{\left(t_0-\frac{1}{n}\right)^{\frac{1}{m_0-1}}} t^{\frac{m_0}{m_0-1}}\, \int_{\Omega}\n(u(t,x))\psi_2(x)\dx \,,
\end{split}
\]
where we have used inequality  \eqref{CP.Monotonicity}\,, in the form  $\n(u(\tau,x))\le \left(t/\tau\right)^{\frac{m_0}{m_0-1}}\n(u(t,x))$ for all $t\ge t_1 +\frac{1}{n}\ge\tau$\,, and the assumptions on $\psi_{1,n}$. Let  $n\to \infty$ to get
\[
\int_{t_0}^{t_1}\int_{\Omega}\n(u(\tau,x))\psi_2(x)\dx \rd\tau \le \frac{m_0-1}{t_0^{\frac{1}{m_0-1}}} t^{\frac{m_0}{m_0-1}}\, \int_{\Omega}\n(u(t,x))\psi_2(x)\dx \qquad\mbox{for all $t\ge t_1 \ge t_0$}\,,
\]
since  $\n(u) \in \LL^1\left((0,\infty):\LL^1_{\Phi_1}(\Omega)\right)$ and  $\|\psi_{1,n}\|_{\LL^\infty(0,+\infty)}\to \|\chi_{[t_0,t_1]}\|_{\LL^\infty(0,+\infty)}=1$\,.

\noindent\textit{Lower estimates. }Let $n$ be so big that $0\le t_0-1/n \le t_1+1/n $. Since  $u\ge 0$\,, we have
\[
\begin{split}
\int_{0}^{\infty}\psi_{1,n}(\tau)\int_\Omega \n(u(\tau,x))\psi_2(x)\dx \rd\tau
&\ge t_0^{\frac{m_0}{m_0-1}}\,\int_{t_0}^{t_1}\frac{\rd\tau}{\tau^{\frac{m_0}{m_0-1}}}\,\int_\Omega \n((t_0,x))\psi_2(x)\dx \\
&=  (m_0-1){t_0^{\frac{m_0}{m_0-1}}} \left[\frac{1}{t_0^{\frac{1}{m_0-1}}}-\frac{1}{t_1^{\frac{1}{m_0-1}}}\right]\,\int_\Omega \n((t_0,x))\psi_2(x)\dx
\end{split}
\]
where we have used used inequality  \eqref{CP.Monotonicity}\,, in the form  $\n(u(t_0-1/n,x))\le \left(\tau/(t_0-1/n)\right)^{\frac{m_0}{m_0-1}}\n(u(\tau,x))$ for all $\tau \ge t_0 -\frac{1}{n}\ge0$\,, together with the fact that $\psi_{1,n}\ge 0$ and $\psi_{1,n}\ge \chi_{[t_0,t_1]}$\,. Since the function $f(t)=t^{-\alpha}$ is convex for  $\alpha=1/(m_0-1)>0$, we may use the inequality $f(t)-f(t_0)\le f'(t)(t-t_0)$ to obtain
\[
\int_{0}^{\infty}\psi_{1,n}(\tau)\int_\Omega \n(u(\tau,x))\psi_2(x)\dx\rd\tau
\ge \left(\frac{t_0}{t_1}\right)^{\frac{m_0}{m_0-1}}(t_1-t_0)\int_\Omega \n(u(t_0,x))\psi_2(x)\dx\,.
\]
Letting now $n\to \infty$ we obtain
\[
\int_{t_0}^{t_1}\int_\Omega \n(u(\tau,x))\psi_2(x)\dx\rd\tau
\ge \left(\frac{t_0}{t_1}\right)^{\frac{m_0}{m_0-1}}(t_1-t_0)\int_\Omega \n(u(t_0,x))\psi_2(x)\dx\,,
\]
which is justified as before. The claim is proved.

\medskip

Summing up the results of the first  steps, \textsl{for every $0\le t_0\le t_1 \le t$\,, and for all $\psi_2(x)$\,, with $\psi_2: \overline{\Omega}\to \RR$ measurable and $\psi_2/\Phi_1$ bounded,  we have}
\begin{equation}\label{step.2.2.thm.repr}\begin{split}
\left(\frac{t_0}{t_1}\right)^{\frac{m_0}{m_0-1}}(t_1-t_0)\int_\Omega \n(u(t_0,x))\psi_2(x)\dx
&\le\int_{\Omega}u(t_0,x)\AI \psi_2(x)\dx - \int_{\Omega}u({t_1},x)\AI \psi_2(x)\dx \\
&\le \frac{m_0-1}{t_0^{\frac{1}{m_0-1}}} t^{\frac{m_0}{m_0-1}}\,\int_\Omega \n(u(t,x))\psi_2(x)\dx\,.
\end{split}\end{equation}

\noindent$\bullet~$\noindent\textsc{Step 4. } We will now prove formula \eqref{thm.NLE.PME.estim} by approximating the kernel of $\AI$,  $\K(x_0,\cdot)$ by means of a sequence of admissible test functions $\psi_{2,n}^{(x_0)}$.   For fixed $x_0\in\Omega$ we consider a sequence of test functions $\psi_{2,n}^{(x_0)}$ with $\psi_{2,n}^{(x_0)}: \overline{\Omega}\to \RR$ measurable and such that $\psi_{2,n}^{(x_0)}/\p$ is bounded\,, such that $\psi_{2,n}^{(x_0)}\to \delta_{x_0}$ as $n\to \infty$\,, in the sense of Radon measures. More specifically, we can choose $\psi_{2,n}^{(x_0)}(x)= |B_{1/n}(x_0)|^{-1}\,\chi_{B_{1/n}(x_0)}(x)$. Therefore, $\AI \psi_{2,n}^{(x_0)}\to \K(\cdot , x_0)$ in $\LL^q(\Omega)$ for all $0<q<N/(N-2s)$.   Then, for any fixed $\tau\ge 0$ we have:
\begin{equation}\label{step.3.1.thm.repr}
\lim_{n\to\infty}  \int_{\Omega}\n(u(\tau,x))\psi_{2,n}^{(x_0)}(x)\dx
= \lim_{n\to\infty} |B_{1/n}(x_0)|^{-1}  \int_{B_{1/n}(x_0)}\n(u(\tau,x))\dx
=   \n(u(\tau,x_0))
\end{equation}
if $x_0$ is a Lebesgue point of the function $x\mapsto u(\tau,x)$\,; $\n(u(\tau,x_0))$ is the corresponding Lebesgue value. If we apply this limit process at the points $\tau=t_0$ and $\tau=t_1$ we get for almost every $x_0$
\begin{equation*}\begin{split}
0\le \left(\frac{t_0}{t_1}\right)^{\frac{m_0}{m_0-1}} (t_1-t_0)\n(u(t_0,x_0))
&\le\lim_{n\to\infty}\int_{\Omega}u(t_0,x)\AI \psi_{2,n}^{(x_0)}(x)\dx - \int_{\Omega}u({t_1},x)\AI \psi_{2,n}^{(x_0)}(x)\dx \\
&\le \frac{m_0-1}{t_0^{\frac{1}{m_0-1}}} t_1^{\frac{m_0}{m_0-1}}\,\n(u(t_1,x_0))<+\infty\,,
\end{split}\end{equation*}
Finally,  since $u\in \mathcal{S}_p$\,, then $u(t)\in \LL^p(\Omega)$ for all $t> 0$\,, with $p>N/(2s)$, and we have already seen in Step 2 that we have
\[
\int_{\Omega}u(t,x)\AI \psi_{2,n}^{(x_0)}(x)\dx \xrightarrow[n\to\, \infty]{} \int_{\Omega}u(t,x)\K(x,x_0)\dx\,,
\]
and formula \eqref{thm.NLE.PME.estim} follows for $t=t_1$\,. For $t$ larger than $t_1$ we use again monotonicity inequality  \eqref{CP.Monotonicity}.\qed


\subsubsection{Consequences of (N1) assumptions}\label{App.F.N1}
Let us first remark that hypothesis (N1) can be rewritten in the following equivalent forms:
\begin{enumerate}
\item[(N1)] $\n\in C^1(\RR\setminus\{0\})$ and $\n/\n'\in {\rm Lip}(\RR)$ and there exists $\mu_0,\mu_1>0$ such that
\[
1-\mu_1\le \left(\frac{\n}{\n'}\right)'\le 1-\mu_0\,,
\]
where $\n/\n'$ is understood to vanish if $\n(r)=\n'(r)=0$ or $r=0$\,.
\item[(N2)] $\n\in C^1(\RR\setminus\{0\})$ and $\n'\in {\rm Lip}_{\rm loc}(\RR\setminus\{0\})$ and there exists $\mu_0,\mu_1>0$ such that
\begin{equation}\label{hyp.1.nonlinearity}
\mu_0\le \frac{\n(r)\n''(r)}{[\n'(r)]^2}\le \mu_1 \qquad\mbox{a.e. }r>0\,.
\end{equation}
\item[(N3)] The function $r\mapsto [\n(r)]^{1-\mu_0}$ is convex ($\log|\n(r)|$\,, if $\mu_0=1$)\,, and $r\mapsto [\n(r)]^{1-\mu_1}$ is concave ($\log|\n(r)|$\,, if $\mu_1=1$) on each of $(-\infty,0)$ and $(0,\infty)$\,.
\end{enumerate}

\begin{lem}[Consequences of hypothesis (N1)]\label{Lem.N1.F}
If a function $H:\RR\to \RR$ satisfies (N1) (in the equivalent form (N2)) with   $0<\alpha_0\le \alpha_1<1$,  i.e. $H\in C^1(\RR\setminus\{0\})$ and $H'\in {\rm Lip}_{\rm loc}(\RR\setminus\{0\})$ and there exist $0<\alpha_0\le\alpha_1<1$ such that
\begin{equation}\label{hyp.1.nonlinearity.H}
\alpha_0\le \frac{H(r)H''(r)}{[H'(r)]^2}\le \alpha_1 \qquad\mbox{a.e. }r>0\,.
\end{equation}
Then the following estimates hold true:
\begin{equation}\label{Lem.N1.1}
\frac{r^{a_0}}{r_0^{a_0}}\le \frac{H(r)}{H(r_0)}\le \frac{r^{a_1}}{r_0^{a_1}}\qquad\mbox{for all $r\ge r_0\ge 0$, with }a_i=\frac{1}{1-\alpha_i}
\end{equation}
and
\begin{equation}\label{Lem.N1.2}
\kb \frac{r^{a_1}}{r_0^{a_1}}\le \frac{H(r)}{H(r_0)}
\le \ka \frac{r^{a_0}}{r_0^{a_0}} \qquad\mbox{for all $0\le r\le r_0$, with }a_i=\frac{1}{1-\alpha_i}
\end{equation}
where $\ka,\kb>0$ and
\begin{equation}\label{Lem.N1.3}
\ka=\left(\frac{1-\alpha_0}{1-\alpha_1}\right)^{a_0}=\left(\frac{a_1}{a_0}\right)^{a_0}>1\,
\qquad\mbox{and}\qquad \kb=\left(\frac{1-\alpha_1}{1-\alpha_0}\right)^{a_1}=\left(\frac{a_0}{a_1}\right)^{a_1}<1\,\,.
\end{equation}
\end{lem}
\noindent {\sl Proof.~}Hypothesis (N1), i.e. inequality \eqref{hyp.1.nonlinearity.H} implies that for a.e. $r>0$
\begin{equation}\label{hyp.1.nonlinearity.H.1}
\alpha_0\frac{H'(r)}{H(r)}\le \frac{H''(r)}{H'(r)}\le \alpha_1\frac{H'(r)}{H(r)}\,,
\end{equation}
an integration on $[r_0,r]$ gives
\begin{equation}\label{hyp.1.nonlinearity.H.2}
\frac{H(r)^{\alpha_0}}{H(r_0)^{\alpha_0}}\le \frac{H'(r)}{H'(r_0)}\le \frac{H(r)^{\alpha_1}}{H(r_0)^{\alpha_1}}
\end{equation}
which is equivalent to the following monotonicity of $H'/H^{\alpha_i}$: if $r_0\le r$ we have
\begin{equation}\label{Lem.N1.H.1}
\frac{H'(r_0)}{H(r_0)^{\alpha_0}}\le \frac{H'(r)}{H(r)^{\alpha_0}}\qquad\qquad\mbox{and}\qquad\qquad \frac{H'(r_0)}{H(r_0)^{\alpha_1}}\ge \frac{H'(r)}{H(r)^{\alpha_1}}\,.
\end{equation}
Integrating both inequalities in two cases, namely when $r\le r_0$and when $r\ge r_0$ will give the desired inequalities \eqref{Lem.N1.2} and \eqref{Lem.N1.2} as follows. Before doing that, we recall that hypothesis (N1) is equivalent to (N3) and implies that
\begin{equation}\label{hyp.1.nonlinearity.H.N3}
(1-\alpha_1)r\le \frac{H(r)}{H'(r)}\le (1-\alpha_0)r\qquad\mbox{or}\qquad
\frac{1}{(1-\alpha_0)r}\le \frac{H'(r)}{H(r)}\le \frac{1}{(1-\alpha_1)r}
\end{equation}

\noindent$\bullet~$\textit{Case I. $r\ge r_0$. Proof of $\eqref{Lem.N1.1}$. }We integrate the first inequality of \eqref{Lem.N1.H.1} on $[r_0,r]$ to get:
\[
\frac{H^{1-\alpha_0}(r)-H^{1-\alpha_0}(r_0)}{1-\alpha_0}\ge (r-r_0)\frac{H'(r_0)}{H(r_0)^{\alpha_0}}
\ge (r-r_0)\frac{H'(r_0)}{H(r_0)}H(r_0)^{1-\alpha_0}\ge \frac{r-r_0}{(1-\alpha_0)r_0}H(r_0)^{1-\alpha_0}
\]
where in the last step we have used \eqref{hyp.1.nonlinearity.H.N3}. It follows that
\begin{equation}\label{Lem.N1.H.2}
H^{1-\alpha_0}(r)\ge H^{1-\alpha_0}(r_0)+\frac{r-r_0}{r_0}H(r_0)^{1-\alpha_0}=H(r_0)^{1-\alpha_0}\frac{r}{r_0}
\end{equation}
that is $H(r)\ge H(r_0)(r/r_0)^{1/(1-\alpha_0)}$ when $r\ge r_0$.\\
On the other hand, integrating the second inequality of \eqref{Lem.N1.H.1} on $[r_0,r]$ we obtain
\[
\frac{H^{1-\alpha_1}(r)-H^{1-\alpha_1}(r_0)}{1-\alpha_1}\le (r-r_0)\frac{H'(r_0)}{H(r_0)^{\alpha_1}}
\le (r-r_0)\frac{H'(r_0)}{H(r_0)}H(r_0)^{1-\alpha_1}\le \frac{r-r_0}{(1-\alpha_1)r_0}H(r_0)^{1-\alpha_1}
\]
where in the last step we have used \eqref{hyp.1.nonlinearity.H.N3}. It follows that
\begin{equation}\label{Lem.N1.H.3}
H^{1-\alpha_1}(r)\le H^{1-\alpha_1}(r_0)+\frac{r-r_0}{r_0}H(r_0)^{1-\alpha_1}=H(r_0)^{1-\alpha_1}\frac{r}{r_0}
\end{equation}
that is $H(r)\le H(r_0)(r/r_0)^{1/(1-\alpha_1)}$ when $r\ge r_0$. Summing up, we have proved \eqref{Lem.N1.1}\,.

\noindent$\bullet~$\textit{Case II. $0\le r\le r_0$. Proof of $\eqref{Lem.N1.2}$. }We integrate the first inequality of \eqref{Lem.N1.H.1} on $[0,r]$ to get
\[
\frac{H^{1-\alpha_0}(r)}{1-\alpha_0}\le \frac{H'(r_0)}{H(r_0)^{\alpha_0}}r
= \frac{H'(r_0)}{H(r_0)}H(r_0)^{1-\alpha_0}r\le  \frac{H(r_0)^{1-\alpha_0}}{1-\alpha_1}\,\frac{r}{r_0}
\]
where in the last step we have used \eqref{hyp.1.nonlinearity.H.N3}. It follows that for all $0\le r\le r_0$ we have
\begin{equation}\label{Lem.N1.H.4}
H(r)\le \left(\frac{1-\alpha_0}{1-\alpha_1}\,\frac{r}{r_0}\right)^{\frac{1}{1-\alpha_0}} H(r_0)\,.
\end{equation}
On the other hand, integrating the second inequality of \eqref{Lem.N1.H.1} on $[0,r]$ to get
\[
\frac{H^{1-\alpha_1}(r)}{1-\alpha_1}\ge \frac{H'(r_0)}{H(r_0)^{\alpha_1}}r
= \frac{H'(r_0)}{H(r_0)}H(r_0)^{1-\alpha_1}r\ge  \frac{H(r_0)^{1-\alpha_1}}{1-\alpha_0}\,\frac{r}{r_0}
\]
where in the last step we have used \eqref{hyp.1.nonlinearity.H.N3}. It follows that for all $0\le r\le r_0$ we have
\begin{equation}\label{Lem.N1.H.5}
H(r)\ge \left(\frac{1-\alpha_1}{1-\alpha_0}\,\frac{r}{r_0}\right)^{\frac{1}{1-\alpha_1}} H(r_0)\,.
\end{equation}
Summing up, we have proved \eqref{Lem.N1.2}\,, and the proof is complete\,.\qed

\section{Comments and extensions}

\noindent$\bullet$ \textbf{Positivity, boundary behaviour and Harnack inequalities. }These are more advanced topics in the  spirit of this paper that will be covered in an upcoming publication \cite{BV-ppr2-2}\,.

\noindent $\bullet$ \textbf{Solutions with any sign. }The abstract functional theory allows to construct solutions for general data not necessarily nonnegative. Due to the property of comparison, extending the upper bounds  to solutions of any sign is easy. Here is the argument: recall first that if $u$ is a solution\,, also $-u$ is a solution. Then, consider the nonnegative solution $u^+$ corresponding to $u_0^+=\max\{u_0,0\}$. Thus by comparison, $u\le u^+$\,, since $u_0\le u_0^+$. Consider also the nonnegative solution $u_-$ corresponding to $u_0^-=-\min\{u_0,0\}$. Then by comparison, $-u\le u^-$\,, since $-u_0\le -u_0^-$, and we get  $-u^-\le u\le u^+$, so that $|u|\le \max\{u^+\,,\,u^-\}$\,.\normalcolor
\normalcolor
We are not considering the detailed theory of signed solutions in this paper because our technique uses in a strong way the monotonicity property \eqref{CP.Monotonicity}, that is only available for nonnegative solutions.

\noindent$\bullet$ \noindent\textbf{RFL and SFL on Lipschitz and non smooth domains. }We have always assumed that the domain is smooth, namely $C^2$, or at least $C^{1,1}$. In the case of Lipschitz domains, for the case of both RFL and SFL there exist (K4) estimates of the Green functions, cf. \cite{Bogdan, Davies1} but $\gamma\in [0,1]$ may vary from the case of regular domains, the constants of the bounds depend strongly on the domain, and the distance to the boundary has to be interpreted in a suitable way. We do not enter this issue in this paper, but we point out that the method presented here could be used to treat the case of non-smooth domains.

\noindent$\bullet$  Let us also mention the possible extension to unbounded domains, and to equations on Riemannian manifolds.

\noindent$\bullet$  \textbf{Other extensions.}  Apart from concentrating on a more detailed analysis of some of the examples that we have proposed in Section \ref{Sect.SFL+RFL+Examples}, there is also the natural extension to nonlinearities of fast diffusion type. As in the standard Laplacian case, we expect the estimates to be considerably different, leading to phenomena like  finite-time extinction, see \cite{BGV-Domains, BV, BV-ADV, BV2012, JLVSmoothing}.

Another possible extension consists in replacing the power function $\lambda^s$  with more general Bernstein functions, as Prof. R. Song has suggested to us. Indeed, paper \cite{Song-subordinate2} contains similar Green function estimates for a class of Bernstein functions. so that assumption $(K2)$ holds in a slightly modified form.

\medskip

\noindent {\textbf{\large \sc Acknowledgment.} Both authors have been partially funded by Project MTM2011-24696 and MTM2014-52240-P (Spain). The authors want to express their gratitude to Professor R. Song for helpful comments about this manuscript.

\vskip .2cm



\addcontentsline{toc}{section}{~~~References}
\begin{thebibliography}{00}
\small

\bibitem{AC} I. Athanasopoulos, L. A. Caffarelli. \textit{Continuity of the temperature in boundary heat control problems} {\rm Adv. Math.} 224 (2010), no. 1, 293-315.

\bibitem{BCr}{\rm  P. B\'enilan, M.~G. Crandall.} \textit{Regularizing effects of homogeneous evolution equations,} {\rm
Contributions to Analysis and Geometry} (suppl.  to Amer. Jour. Math.), Johns Hopkins Univ. Press, Baltimore, Md., 1981. Pp.
23-39.

\bibitem{BlGe} R.~M. Blumenthal, R.~K. Getoor. \emph{Some theorems on stable processes.} Trans. Amer. Math. Soc. 95 (1960), no.~2, 263--273.

\bibitem{bogdan-censor} {\rm K. Bogdan, K. Burdzy, K.,  Z.-Q. Chen}. \emph{ Censored stable processes}.  Probab. Theory Relat. Fields {\bf 127} (2003), 89--152.

\bibitem{Bogdan} K. Bogdan, T. Grzywny, M. Ryznar\textit{Heat Kernel Estimates for the Fractional Laplacian with Dirichlet conditions.} Ann. Probab. \textbf{38} 2010, No. 5, 1901--1923

\bibitem{BGV-Domains} {\rm M. Bonforte, G. Grillo, J. L. V{\'a}zquez.} {\it Behaviour near extinction for the Fast Diffusion Equation on bounded domains}, \textrm{J. Math. Pures Appl.} \textbf{97} (2012), 1-–38.\normalcolor

\bibitem{BV} {\rm M.  Bonforte, J.~L. V\'azquez}. \textit{Global positivity estimates and Harnack inequalities for the fast diffusion equation,} \textrm{J. Funct. Anal.} \textbf{240}, no.2, (2006) 399--428.

\bibitem{BSV2013} M. Bonforte, Y. Sire, J.~L. V\'azquez, \textit{Existence, Uniqueness and Asymptotic behaviour for fractional porous medium on bounded domains.} Discr. Cont. Dyn. Sys. \textbf{35} (2015), 5725--5767.

\bibitem{BV-ADV}  {\rm M.  Bonforte, J.~L. V\'azquez. } {\it Positivity, local smoothing, and Harnack inequalities for very fast diffusion equations}, \textrm{Advances in Math.} \bf 223 \rm (2010), 529--578.

\bibitem{BV2012} M. Bonforte, J.~L. V\'azquez. \textit{Quantitative Local and Global  A Priori Estimates for Fractional Nonlinear Diffusion Equations}, Advances in Math \textbf{250} (2014), 242--284. 

\bibitem{BV-PPR1} M. Bonforte,  J. L. V\'azquez, \emph{A Priori Estimates  for Fractional Nonlinear  Degenerate Diffusion Equations on bounded domains},  Arch. Rat. Mech. Anal. \textbf{218} (2015), no. 1, 317--362.

\bibitem{BV-ppr2-2} M. Bonforte, J.~L. V\'azquez, \textit{Fractional Nonlinear Degenerate Diffusion Equations on Bounded Domains. Part II}, in preparation.

\bibitem{Brezis71} H. Brezis.\textit{Monotonicity methods in Hilbert spaces and some applications to nonlinear partial differential equations}, Proc. Symp. Nonlinear Funct. Anal., Madison, Acad. Press (1971), 101--156.

\bibitem{Cabre-Tan}X. Cabr\'e, J. Tan. {\it Positive solutions of nonlinear problems involving the square root of the Laplacian}, Adv. in Math. \textbf{224} (2010), no.~5, 2052--2093.

\bibitem{Caffarelli-Silvestre} L. Caffarelli, L. Silvestre, \textit{An extension problem related to the fractional Laplacian}.  Comm. Partial Diff. Eq.  \textbf{32} (2007),
no.~7-9, 1245--1260.

\bibitem{Song-intrinsic-ultra}Z.-Q. Chen, R. Song, \textit{Intrinsic ultracontractivity and conditional gauge for symmetric stable processes.} J. Funct. Anal., \textbf{150} (1997), 204--239.

\bibitem{Song-Sum3}Z.-Q. Chen, P. Kim, R. Song, \textit{Dirichlet heat kernel estimates for $\Delta^{\alpha/2}+\Delta^{\beta/2}$}. Illinois J. Math. \textbf{54} (2010), no. 4, 1357--1392.

\bibitem{Song-coeff}Z.-Q. Chen, P. Kim, R. Song, \textit{Two-sided heat kernel estimates for censored stable-like processes. }Probab. Theory Relat. Fields \textbf{146} (2010) 361--399

\bibitem{Song-Rel}Z.-Q. Chen, P. Kim, R. Song, \textit{Green function estimates for relativistic stable processes in half-space-like open sets.} Stochastic Process. Appl. \textbf{121} (2011), no. 5, 1148--1172.

\bibitem{Song-Drift}Z.-Q. Chen, P. Kim, R. Song, \textit{Dirichlet heat kernel estimates for fractional Laplacian with gradient perturbation. }Ann. Probab. \textbf{40} (2012), no. 6, 2483-–2538.

\bibitem{Song-Sum2}Z.-Q. Chen, P. Kim, R. Song, Z. Vondra{\v{c}}ek. \textit{Sharp Green function estimates for $\Delta+\Delta^{\alpha/2}$ in $C^{1,1}$ open sets and their applications.} Illinois J. Math. \textbf{54} (2010), no. 3, 981--1024.

\bibitem{CL71} { M.~G. Crandall, T.M. Liggett.} \textit{Generation of semi-groups of nonlinear transformations on general
Banach spaces},  {Amer. J. Math.} {\bf 93} (1971) 265--298.

\bibitem{CP-JFA} M.~G. Crandall, M. Pierre, \textit{Regularizing Effectd for $u_t=A\varphi(u)$ in $\LL^1$}, J. Funct. Anal. \textbf{45}, (1982), 194-212

\bibitem{DK0} {\rm B. Dahlberg, C.~E. Kenig}. \textit{Nonnegative solutions of the initial-Dirichlet problem for generalized porous medium equation in cylinders}, \textrm{J. Amer. Math. Soc.} \textbf{1} (1988), 401-–412.

\bibitem{DaskaBook}P. Daskalopoulos, C. E. Kenig, \textsl{``Degenerate diffusions. Initial value problems and local regularity theory''. }EMS Tracts in Mathematics, \textbf{1}. European Mathematical Society (EMS), Zürich, 2007. x+198 pp. ISBN: 978-3-03719-033-3

\bibitem{D1} E. B. Davies.  \textit{Explicit constants for Gaussian upper bounds on heat kernels}, Amer. J. Math. \textbf{109} (1987), no. 2, 319–333.

\bibitem{D2} E. B. Davies.  \textit{The equivalence of certain heat kernel and Green function bounds}, J. Funct. Anal. \textbf{71} (1987), no. 1, 88–103.

\bibitem{Davies1}{\rm E. B. Davies}. {\sl ``Heat kernels and spectral theory'',} Cambridge Tracts in Mathematics, 92. Cambridge University Press, Cambridge, 1990. x+197 pp. ISBN: 0-521-40997-7

\bibitem{Davies2} {\rm E. B. Davies}. {\sl ``Spectral theory and differential operators'',} Cambridge Studies in Advanced Mathematics, 42. Cambridge University Press, Cambridge, 1995. x+182 pp. ISBN: 0-521-47250-4

\bibitem{DS} E. B. Davies, B. Simon. \textit{Ultracontractivity and the heat kernel for Schrödinger operators and Dirichlet Laplacians}, J. Funct. Anal. \textbf{59} (1984), no. 2, 335–395.

\bibitem{Jak} T. Jakubowski. \textit{The estimates for the Green function in Lipschitz domains for the symmetric stable processes. }
Probab. Math. Statist. \textbf{22} (2002), no. 2, Acta Univ. Wratislav. No. 2470, 419--441.

\bibitem{Kim-Coeff} K.-Y. Kim, P. Kim. \textit{Two-sided estimates for the transition densities of symmetric Markov processes dominated by stable-like processes in $C^{1,\eta}$ open sets}, Stochastic Processes and their Applications {\bf 124} (2014) 3055--3083.

\bibitem{Song-subordinate1} P. Kim, H. Park, R. Song. \textit{Sharp estimates on the Green functions of perturbations of subordinate Brownian motions in bounded $\kappa$-fat open sets. }Potential Analysis, \textbf{38} (2013), 319--344

\bibitem{song-NonSymm-1}P. Kim, R. Song. \textit{Estimates on Green functions and Schrödinger-type equations for non-symmetric diffusions with measure-valued drifts.} J. Math. Anal. Appl. \textbf{332} (2007), no. 1, 57--80.

\bibitem{Song-Sum1Gen}P. Kim, R. Song, Z. Vondra{\v{c}}ek. \textit{Potential theory of subordinate Brownian motions with Gaussian components. } Stochastic Process. Appl. \textbf{123} (2013), 764--795.

\bibitem{Song-subordinate2} P. Kim, R. Song, Z. Vondra{\v{c}}ek. \textit{Two-sided Green function estimates for killed subordinate Brownian motions.} Proc. Lond. Math. Soc. (3) \textbf{104} (2012), no. 5, 927--958.

\bibitem{Kul} T. Kulczycki. \textit{Properties of Green function of symmetric stable processes. }Probab. Math. Statist. \textbf{17} (1997), no. 2, Acta Univ. Wratislav. No. 2029, 339--364.

\bibitem{DPQRV1} {\rm A.~de Pablo, F. Quir\'os, A. Rodr\'iguez, J. L. V\'azquez. }\textit{A fractional porous medium equation} Adv. Math. \textbf{226} (2011), no. 2, 1378–1409.

\bibitem{DPQRV2}{\rm A.~de Pablo, F. Quir\'os, A. Rodr\'iguez, J. L. V\'azquez. } \textit{A general fractional porous medium equation},  Comm. Pure Applied Math. {\bf 65} (2012), 1242--1284.

\bibitem{Pierre} {\rm  M. Pierre}.   \textit{Uniqueness of the solutions of $u_t-\Delta \phi(u)=0$ with initial datum a measure}. \textrm{ Nonlinear Anal. T. M. A.}  {\bf 6} (1982), pp. 175--187.

\bibitem{RosSer} X. Ros-Oton, J. Serra, \textit{The Dirichlet problem for the fractional Laplacian: regularity up to the boundary}, J. Math. Pures Appl. \textbf{101} (2014), 275--302.

\bibitem{JLVmonats}{\rm  J.~L. V{\'a}zquez. } \textit{The Dirichlet Problem for the Porous Medium Equation in Bounded Domains. Asymptotic Behavior}, Monatsh. Math. \textbf{142}, (2004) 81–111 .

\bibitem{VazBook}{\rm  J.~L. V{\'a}zquez. } {\sl ``The Porous Medium Equation. Mathematical Theory''}, vol.~Oxford Mathematical Monographs, Oxford University Press, Oxford, 2007.

\bibitem{JLVSmoothing}{\rm J.~L. V\'azquez}. {\sl ``Smoothing and Decay Estimates for Nonlinear Diffusion Equations. Equations of Porous Medium
Type''.} Oxford Lecture Series in Mathematics and its Applications, 33. Oxford University Press, Oxford, 2006.

\bibitem{Vaz2012}{\rm  J.~L. V{\'a}zquez. } \textit{Barenblatt solutions and asymptotic behaviour for a nonlinear fractional heat equation of porous medium type}, J. Eur. Math. Soc. {\bf 16} (2014), 769--803. 

\bibitem{Zh2002} Q. S. Zhang.  \textit{The boundary behavior of heat kernels of Dirichlet Laplacians}, J. Differential Equations \textbf{182} (2002), no. 2, 416–430.


\end{thebibliography}
\end{document}